\documentclass[11pt, a4paper]{amsart}

\usepackage{amsfonts,amsmath,amssymb, amscd}

\title{Polynomial identities and noncommutative versal torsors}

\author{Eli Aljadeff}
\address{Eli Aljadeff:
Department of Mathematics,
Technion - Israel Institute of Technology,
32000 Haifa, Israel}
\email{aljadeff@techunix.technion.ac.il}

\author{Christian Kassel}
\address{Christian Kassel: Institut de Recherche Math\'{e}matique Avanc\'{e}e,
CNRS - Universit\'{e} Louis Pasteur,
7 rue Ren\'{e} Descartes, 67084 Strasbourg, France}
\email{kassel@math.u-strasbg.fr}

\newtheorem{Lem}{Lemma}[section]
\newtheorem{Prop}[Lem]{Proposition}
\newtheorem{Cor}[Lem]{Corollary}
\newtheorem{Thm}[Lem]{Theorem}

\theoremstyle{definition}

\newtheorem{Rem}[Lem]{Remark}
\newtheorem{Rems}[Lem]{Remarks}
\newtheorem{Expl}[Lem]{Example}

\newcommand\pf{\begin{proof}}
\newcommand\epf{\end{proof}}

\newcommand\Alg{\operatorname{Alg}}

\newcommand\End{\operatorname{End}}
\newcommand\Hom{\operatorname{Hom}}
\newcommand\Tr{\operatorname{Tr}}
\newcommand\Ker{\operatorname{Ker}}

\newcommand\Forms{\operatorname{Forms}}

\newcommand\id{{\operatorname{id}}}

\renewcommand\o{\otimes}

\newcommand\RComod[1]{\text{Hom}^{#1}}

\newcommand\RAlgComod[1]{\Alg^{#1}}

\renewcommand\AA{\mathcal A}
\newcommand\BB{\mathcal B}
\newcommand\UU{\mathcal U}
\newcommand\ZZ{\mathcal Z}
\newcommand\KK{\mathcal K}

\newcommand\mm{\mathfrak{m}}
\newcommand\MM{\mathfrak{M}}

\newcommand\eps{\varepsilon}

\newcommand\sw[1]{{}_{(#1)}}

\numberwithin{equation}{section}

\begin{document}

\maketitle

\vskip 30pt
\noindent
{\sc Abstract}. 
\emph{ 
To any cleft Hopf Galois object, i.e., any algebra~${}^{\alpha} H$
obtained from a Hopf algebra~$H$
by twisting its multiplication with a two-cocycle~$\alpha$,
we attach two ``universal algebras''
$\AA_H^{\alpha}$ and~$\UU_H^{\alpha}$.
The algebra $\AA_H^{\alpha}$ is obtained by twisting the multiplication
of~$H$ with the most general two-cocycle~$\sigma$ formally
cohomologous to~$\alpha$. 
The cocycle~$\sigma$ takes values in the field of rational functions on~$H$.
By construction, $\AA_H^{\alpha}$ is a cleft $H$-Galois extension of
a ``big'' commutative algebra~$\BB_H^{\alpha}$.
Any ``form'' of~${}^{\alpha} H$ can be obtained 
from~$\AA_H^{\alpha}$ by a specialization of~$\BB_H^{\alpha}$
and vice versa. If the algebra ${}^{\alpha} H$ is simple, then
$\AA_H^{\alpha}$ is an Azumaya algebra with center~$\BB_H^{\alpha}$.
The algebra~$\UU_H^{\alpha}$ is constructed using a general
theory of polynomial identities that we set up for arbitrary comodule algebras;
it is the universal comodule algebra 
in which all comodule algebra identities of~${}^{\alpha} H$ are satisfied.
We construct an embedding of $\UU_H^{\alpha}$ into~$\AA_H^{\alpha}$;
this embedding maps the center~$\ZZ_H^{\alpha}$ of $\UU_H^{\alpha}$
into~$\BB_H^{\alpha}$ when the algebra ${}^{\alpha} H$ is simple.
In this case, under an additional assumption, 
$\AA_H^{\alpha} \cong \BB_H^{\alpha} \otimes_{\ZZ_H^{\alpha}} \UU_H^{\alpha}$,
thus turning~$\AA_H^{\alpha}$
into a central localization of~$\UU_H^{\alpha}$.
We completely work out these constructions in the case of the four-dimensional
Sweedler algebra.
}

\bigskip
\noindent
{\sc Key Words:}
Hopf algebra, comodule algebra, Galois extension, cocycle, polynomial identity

\bigskip
\noindent
{\sc Mathematics Subject Classification (2000):}
16R50, 16W30, 16S35, 16S38, 16S40, 16H05, 16E99, 17B37, 55R10, 58B32, 81R50, 81R60

\bigskip\bigskip

\hspace{3cm}

\section*{Introduction}

In this paper we deal with algebras~${}^{\alpha} H$
obtained from a Hopf algebra~$H$
by twisting its multiplication with a two-cocycle~$\alpha$. 
This class of algebras coincides with the class of so-called 
\emph{cleft Hopf Galois extensions} of the ground field; 
strongly $G$-graded algebras and certain $G$-algebras, where $G$ is a finite group,
belong to this class. 
As has been stressed many times (see e.g.~\cite{S}), Hopf Galois
extensions can be viewed as noncommutative analogues of $G$-torsors
or of principal fiber bundles 
where the role of the structural group is played by a Hopf algebra. 
Hopf Galois extensions abound in the world of quantum groups
and of noncommutative geometry
(see, e.g., \cite{BDZ}, \cite{BM1}, \cite{BM2}, 
\cite{Ha1}, \cite{HM}, \cite{LPR}).

To the algebra~${}^{\alpha} H$ we attach two ``universal algebras''
$\UU_H^{\alpha}$ and~$\AA_H^{\alpha}$, 
using two constructions of a very different nature.
The algebra~$\UU_H^{\alpha}$ is a graded quotient of the tensor algebra
over the underlying vector space of~$H$.
To construct this quotient, we set up a theory of 
\emph{polynomial identities for comodule algebras}.
We define~$\UU_H^{\alpha}$ as the universal comodule algebra 
in which all $H$-comodule algebra identities of~${}^{\alpha} H$ are satisfied.

The second universal algebra,~$\AA_H^{\alpha}$, is obtained by twisting the multiplication
of~$H$ with a two-cocycle~$\sigma$ taking
values in the field of rational functions on~$H$.
The cocycle~$\sigma$ can be viewed as the most general cocycle
formally cohomologous to~$\alpha$.
By construction, $\AA_H^{\alpha}$ is a cleft $H$-Galois extension of
the commutative algebra~$\BB_H^{\alpha}$ generated by the values of
the cocycle~$\sigma$ and of its convolution inverse.
We show that any ``form'' of~${}^{\alpha} H$ can be obtained from~$\AA_H^{\alpha}$
by a specialization of~$\BB_H^{\alpha}$. 
Conversely, any central specialization of~$\AA_H^{\alpha}$
is a form of~${}^{\alpha} H$.
Thus, the set of algebra morphisms~$\Alg(\BB_H^{\alpha},K)$
parametrizes the isomorphism classes of $K$-forms of~${}^{\alpha} H$
and $\AA_H^{\alpha}$ can be viewed as a flat deformation of~${}^{\alpha} H$
over the commutative algebra~$\BB_H^{\alpha}$; it is a noncommutative
analogue of a \emph{versal deformation space} or a \emph{versal torsor} in the sense of Serre
(see~\cite[Chap.~I]{GMS}).
We also prove that, if ${}^{\alpha} H$ is a simple algebra, then
$\AA_H^{\alpha}$ is an Azumaya algebra with center~$\BB_H^{\alpha}$.

We relate the algebras $\UU_H^{\alpha}$ and~$\AA_H^{\alpha}$
by constructing an injective comodule algebra morphism 
$\mu_{\sigma} : \UU_H^{\alpha} \to \AA_H^{\alpha}$.
The morphism $\mu_{\sigma}$ sends 
the center~$\ZZ_H^{\alpha}$ of~$\UU_H^{\alpha}$ to the center~$\BB_H^{\alpha}$
of~$\AA_H^{\alpha}$ when the algebra ${}^{\alpha} H$ is simple. 
Under an additional assumption
we prove that $\mu_{\sigma}$ induces
an isomorphism of $H$-comodule algebras
$\BB_H^{\alpha} \otimes_{\ZZ_H^{\alpha}} \UU_H^{\alpha} \cong \AA_H^{\alpha}$.
We thus obtain the ``versal deformation space''~$\AA_H^{\alpha}$
as a central localization of the universal $H$-comodule algebra~$\UU_H^{\alpha}$.

Most of the results presented here generalize results obtained in the group case
by the first-named author jointly with D.~Haile and M.~Natapov, see~\cite{AHN1}. 
In this important special case, additional results are obtained in~\cite{AHN1}
such as a characterization of the pairs $(G, \alpha)$ 
consisting of a finite group~$G$ and a two-cocycle~$\alpha$
for which $\AA_H^{\alpha}$~is a division algebra,
see also~\cite{AN}.
Since we deal with arbitrary Hopf algebras,
our results, whose proofs are not straightforward extensions 
of those of~\cite{AHN1},
cover situations that can be radically different from the group case.
Our general framework also leads to the emergence of
new interesting questions on Hopf algebras
such as the one posed in Section~\ref{integrality}.

In~\cite{K},~\cite{KS} a project aimed at classifying
Hopf Galois extensions was started. 
The results presented in this paper can be considered 
as first steps towards the construction of universal Hopf Galois extensions.
Although we do not obtain universal extensions (which may well not exist), 
we manage to construct versal deformation spaces. 
We believe that such deformation spaces are of interest 
and that they deserve to be explicited for many Hopf Galois extensions.

The paper is organized as follows. 
In Section~\ref{basics} we recall some well-known facts on comodule algebras.

In Section~\ref{identities} we define the concept of an $H$-identity
for an $H$-comodule algebra and 
in Section~\ref{univ-comod-alg} we associate 
to each $H$-comodule algebra~$A$ the universal algebra~$\UU_H(A)$.

In Section~\ref{universal evaluation} we restrict ourselves to
the case of $H$-comodule algebras~$A$ that are isomorphic 
as comodules to~$B\otimes H$, 
where $B$ is the subalgebra of coinvariants of~$A$,
assumed to be central in~$A$.
For each such $H$-comodule algebra
we construct a comodule algebra morphism
whose kernel turns out to be exactly the ideal of $H$-identities for~$A$.
This allows us to embed the universal algebra~$\UU_H(A)$ 
into a more controllable algebra, 
namely the tensor product of~$A$ by a (commutative) polynomial algebra.

We construct the two-cocycle~$\sigma$ and
the commutative algebra~$\BB_H^{\alpha}$ in Section~\ref{universal cocycle}.
In Section~\ref{universal twisted algebra} we define the ``universal
twisted algebra"~$\AA_H^{\alpha}$ and show that,
if~${}^{\alpha} H$ is (semi)simple, then $\AA_H^{\alpha}$ is ``generically" (semi)simple,
i.e., $\KK_H^{\alpha} \otimes_{\BB_H^{\alpha}} \AA_H^{\alpha}$ is
(semi)simple, where~$\KK_H^{\alpha}$ is the field of fractions of~$\BB_H^{\alpha}$.

In Section~\ref{forms} we define forms of~${}^{\alpha} H$ 
and show how specializations of~$\BB_H^{\alpha}$
give rise to forms of~${}^{\alpha} H$ and \emph{vice versa}. We also prove that
$\AA_H^{\alpha}$ is an Azumaya algebra if ${}^{\alpha} H$ is a simple algebra.

In Section~\ref{forms-identities} we embed 
$\UU_H^{\alpha}$ into~$\AA_H^{\alpha}$
and we show in Section~\ref{nondegenerate} that, if ${}^{\alpha} H$ is simple,
then the previous embedding maps the center~$\ZZ_H^{\alpha}$ 
of~$\UU_H^{\alpha}$ into~$\BB_H^{\alpha}$,
resulting in the isomorphism
$\BB_H^{\alpha} \otimes_{\ZZ_H^{\alpha}} \UU_H^{\alpha}
\cong \AA_H^{\alpha}$.

We illustrate our results in Section~\ref{Sweedler}
on the four-dimen\-sion\-al Sweedler Hopf algebra.
In this example the base of our versal deformation space is
a quadric hypersurface in a five-dimensional affine space
deprived of two hyperplanes.

The paper ends with two appendices
related to the constructions of Section~\ref{universal cocycle}.
In the first appendix we construct a map needed 
for the construction of the cocycle~$\sigma$.
In the second appendix we relate the algebra~$S(t_H)_{\Theta}$
defined in \S~\ref{universal cocycle2}
to Takeuchi's free commutative Hopf algebra over a coalgebra.
It follows that $S(t_H)_{\Theta}$ is obtained from a polynomial algebra
by inverting certain canonical polynomials. These polynomials extend
Dedekind's group determinants to the framework of coalgebras.

\section{Preliminaries}\label{basics}

Throughout the paper, we fix a ground field~$k$ over which all
our constructions are defined. As usual, 
unadorned tensor products mean tensor products over~$k$.

All algebras that we consider are associative unital $k$-algebras.
The unit of an algebra~$A$ will be denoted~$1_A$, or~$1$ if no confusion is possible.
All algebra morphisms are supposed to preserve the units.
We denote the set of algebra morphisms from $A$ to~$A'$ by $\Alg(A,A')$.

All coalgebras are coassociative counital $k$-coalgebras. 
We denote
the comultiplication of a coalgebra by $\Delta$ and its counit by~$\eps$.
We use the Heyneman-Sweedler sigma notation (see \cite[Sect.~2.1.1]{Abe}, \cite[Sect.~1.2]{Sw}):
$$\Delta(x) = \sum_{(x)} x\sw1 \o x\sw2$$
for the comultiplication of $x \in C$ and
$$\Delta^{(2)}(x) = \sum_{(x)} x\sw1 \o x\sw2 \o x_{(3)}$$
for the iterated comultiplication 
$\Delta^{(2)} = (\Delta \o \id_C) \circ \Delta = (\id_C \o \Delta) \circ \Delta$,
and so~on.

Let $H$ be a Hopf algebra
with comultiplication~$\Delta$,
counit~$\eps$, and antipode~$S$.
An \emph{$H$-comodule algebra} is an algebra $A$ equipped with a
right $H$-comodule structure whose (coassociative, counital) coaction
$$\delta : A \to A \o H$$ 
is an algebra morphism
(we will not use left comodules).
A coaction~$\delta$ is called \emph{trivial} if
$\delta(a)  = a \otimes 1$ for all $a \in A$.

Given two $H$-comodule algebras $A$ and $A'$ with respective coactions
$\delta$ and~$\delta'$, an algebra morphism  $f : A \to A'$ 
is an \emph{$H$-comodule algebra morphism} if
$$\delta' \circ f = (f\otimes \id_H) \circ \delta \, .$$
We denote $\RAlgComod{H}$ the category whose objects are $H$-comodule algebras
and arrows are $H$-comodule algebra morphisms.

The space of \emph{co\-\"\i n\-var\-iants} of an $H$-comodule algebra~$A$
is the subspace~$A^H$ of~$A$ defined by
$$A^H = \{ a \in A \, | \, \delta(a)  = a \otimes 1\} \, .$$
The subspace $A^H$ is a subalgebra and a subcomodule of~$A$.
It is the largest subcomodule of~$A$ on which the coaction is trivial.

We shall later use the following lemma.

\begin{Lem}\label{inject-coinv}
Let $f : A_1 \to A_2$ be an injective $H$-comodule algebra morphism.
Then 
$$A_1^H = \{ a \in A_1 \, |\, f(a) \in A_2^H\}\, .$$
\end{Lem}

\pf
If $a \in A_1^H$, then $\delta(a) = a\otimes 1$.
Therefore, $\delta(f(a)) = f(a) \otimes 1$,
which shows that $f(a)$ belongs to~$A_2^H$.

Conversely, let $a\in A_1$ be such that $\delta(f(a)) = f(a) \otimes 1$.
Let $\{1\} \cup \{h_i\}_{i\in I}$ be a basis of~$H$.
Expanding $\delta(a) = a_0 \otimes 1 + \sum_{i\in I} \, a_i \otimes h_i$
in this basis, we obtain
\begin{eqnarray*} 
f(a) \otimes 1 &= & \delta(f(a))
= (f \otimes \id_H) \bigl( \delta(a) \bigr)  \\
&=& (f \otimes \id_H) \Bigl( a_0 \otimes 1 + \sum_{i\in I} \, a_i \otimes h_i \Bigr) \\
&=& f(a_0) \otimes 1 + \sum_{i\in I} \, f(a_i) \otimes h_i\, .
\end{eqnarray*} 
This implies that $f(a_i) = 0$ for all $i\in I$. 
Since $f$ is injective, $a_i = 0$ for all $i\in I$. 
Therefore, $\delta(a) = a_0 \otimes 1$.
Applying $\id_{A_1} \otimes \eps$ to both sides of the previous equality,
we obtain $a = a_0$. It follows that $a \in A_1^H$.
\epf

We end these preliminaries with two well known classes of comodule algebras.

\begin{Expl}\label{graded algebra}
The group algebra $H = k[G]$ of a group~$G$ is a Hopf algebra
with comultiplication, counit, and antipode given for all $g\in G$
by
$$\Delta(g) = g\otimes g \, , \quad\eps(g) = 1\, , \quad S(g) = g^{-1}\, .$$
It is well-known (see~\cite[Lemma~4.8]{BM}) that 
an $H$-comodule algebra~$A$ is the same as a $G$-graded algebra
$$A = \bigoplus_{g\in G}\, A_g \, .$$
The coaction $\delta : A \to A \otimes H$ is given by
$\delta(a) = a\otimes g$ for $a\in A_g$ and~$g\in G$.
We have $A^H = A_e$, where $e$ is the neutral element of~$G$.
\end{Expl}

\begin{Expl}\label{algebra with G-action}
Let $G$ be a finite group and $H = k^G$ be 
the algebra of $k$-valued functions on a finite group~$G$.
This algebra can be equipped with a Hopf algebra structure 
that is dual to the Hopf algebra~$k[G]$ above. 
An $H$-comodule algebra~$A$ is the same as an algebra
equipped with a left action of $G$ on~$A$ by group automorphisms.
If we denote the action of~$g\in G$ on $a\in A$ by ${}^g a$, then
the coaction $\delta : A \to A \otimes H$ is given by
$$\delta(a) = \sum_{g\in G}\, {}^g a \otimes e_g\, ,$$
where $\{e_g\}_{g\in G}$ is the basis of~$H$ consisting of
the functions $e_g$ defined by $e_g(h) = 1$ if $h=g$, and~$0$ otherwise.
The subalgebra of coinvariants of~$A$ coincides with the subalgebra of
$G$-invariant elements: $A^H = A^G$.
\end{Expl}

\section{$H$-identities}\label{identities}

Let $H$ be a Hopf algebra.
We first build up the formal setting in which the identities will live.

\subsection{The tensor algebra}\label{T(XH)}

Let $X_H$ be a copy of the underlying vector space of~$H$.
We denote the identity map from $H$ to~$X_H$ by $h\mapsto X_h$ for all $h\in H$.
By definition we have
$$X_{\lambda_1 h_1 + \lambda_2 h_2} = \lambda_1 X_{h_1} + \lambda_2 X_{h_2}$$
for all $\lambda_1, \lambda_2 \in k$ and $h_1, h_2 \in H$.

The vector space $X_H$ is naturally an $H$-comodule whose coaction
$\delta$ is induced by the comultiplication $\Delta$ of~$H$, i.e.,
for all $h\in H$ we have
\begin{equation}\label{TXcoaction}
\delta(X_h) = \sum_{(h)}\, X_{h_{(1)}} \otimes h_{(2)} \in X_H \otimes H \, .
\end{equation}

Consider the tensor algebra $T(X_H)$ of the vector space~$X_H$
over the ground field~$k$:
$$T(X_H) = \bigoplus_{r\geq 0}\, T^r(X_H) \, ,$$
where $T^r(X_H) = k$ if $r=0$, $T^r(X_H) = X_H$ if $r=1$, and 
$T^r(X_H) = X_H^{\otimes r}$ for all $r\geq 2$. 
If $\{h_i\}_{i\in I}$ is a basis of~$H$, 
then $T(X_H)$ is the free noncommutative algebra
over the set of indeterminates $\{X_{h_i}\}_{i\in I}$.

We equip $T(X_H)$ with an $H$-comodule algebra structure with 
the coaction induced by~\eqref{TXcoaction}.
On~$T^r(X_H)$ ($r \geq 2$) the coaction is given by
\begin{equation}\label{TXcoaction2}
\delta(X_{h^{[1]}} \cdots X_{h^{[r]}})
= \sum_{(h^{[1]}), \ldots, (h^{[r]})}\,
X_{h^{[1]}_{(1)}} \cdots X_{h^{[r]}_{(1)}}
\otimes h^{[1]}_{(2)} \cdots h^{[r]}_{(2)}
\end{equation}
for all $h^{[1]}, \ldots, h^{[r]} \in H$.

By the universal property of $T(X_H)$, any algebra morphism from $T(X_H)$ to an algebra~$A$
is determined by its restriction to its degree-one summand $T^1(X_H) = X_H$.
We thus have a natural bijection
$$\Alg(T(X_H), A) \cong \Hom(X_H,A) = \Hom(H,A)\, .$$
If we restrict this bijection to $H$-comodule algebra morphisms,
we obtain a bijection
$$\RAlgComod{H}(T(X_H), A) \cong \RComod{H}(X_H,A) = \RComod{H}(H,A) \, ,$$
which is functorial in~$A \in \RAlgComod{H}$.
Here $\RAlgComod{H}$ is the category of $H$-comodule algebras
and $\RComod{H}$ is the category of $H$-comodules.

Let us give a typical example of a co\-\"\i n\-var\-i\-ant element of~$T(X_H)$.

\begin{Lem}\label{degree2-identity}
For any $h^{[1]}, \ldots, h^{[r]} \in H$,
$$P_{h^{[1]}, \ldots, h^{[r]}} = \sum_{(h^{[1]}), \ldots, (h^{[r]})}\,
X_{h^{[1]}_{(1)}} \cdots X_{h^{[r]}_{(1)}} X_{S(h^{[1]}_{(2)} \cdots h^{[r]}_{(2)})}$$
is a co\-\"\i n\-var\-i\-ant element of $T(X_H)$.
\end{Lem}

\pf
By \eqref{TXcoaction2}, 
the element $\delta(P_{h^{[1]}, \ldots, h^{[r]}}) \in T(X_H) \otimes H$ is equal to
\begin{multline*}
\sum_{(h^{[1]}), \ldots, (h^{[r]})}\,
X_{h^{[1]}_{(1)}} \cdots X_{h^{[r]}_{(1)}} X_{S(h^{[1]}_{(4)} \cdots h^{[r]}_{(4)})}
\otimes {h^{[1]}_{(2)}} \cdots {h^{[r]}_{(2)}} {S(h^{[1]}_{(3)} \cdots h^{[r]}_{(3)})} \\
= 
\sum_{(h^{[1]}), \ldots, (h^{[r]})}\,
X_{h^{[1]}_{(1)}} \cdots X_{h^{[r]}_{(1)}} X_{S(h^{[1]}_{(3)} \cdots h^{[r]}_{(3)})}
\otimes \eps({h^{[1]}_{(2)}}) \cdots \eps({h^{[r]}_{(2)}}) 1\, ,
\end{multline*}
which is clearly equal to~$P_{h^{[1]}, \ldots, h^{[r]}} \otimes 1$.
\epf

\subsection{$H$-identities for comodule algebras}\label{identities-def}

Let $A$ be an $H$-comodule algebra. An element $P\in T(X_H)$ is called
an \emph{$H$-identity for~$A$} if $\mu(P) = 0$ for all 
$H$-comodule algebra morphisms $\mu : T(X_H) \to A$.

Our $H$-identities for $H$-comodule algebras
generalize graded identities for group-graded algebras.
Indeed, let $H= k[G]$ be the Hopf algebra of a group~$G$.
The set $\{X_g\}_{g\in G}$ forms a basis of~$X_H$.
Let $A$ be an $H$-comodule algebra; 
by~Example~\ref{graded algebra}, it is a
$G$-graded algebra: $A = \bigoplus_{g\in G}\, A_g$. 
An $H$-comodule algebra morphism $\mu : T(X_H) \to A$ is then 
an algebra morphism such that $\mu(X_g) \in A_g$ for all $g\in G$,
and an $H$-identity for~$A$ is what is called a \emph{graded identity}, 
see~\cite{AHN1},~\cite{BZ}. 

Our definition of identities for Hopf comodule algebras
should also be compared to the definition of identities 
for Hopf module algebras appearing in~\cite{BL}.

Let $I_H(A)$ be the set of all $H$-identities for~$A$.
By definition,
\begin{equation}\label{I(A)-def}
I_H(A) = \bigcap_{\mu \, \in \, \RAlgComod{H}(T(X_H),A)}\, \Ker \mu \, .
\end{equation}

\begin{Prop}\label{Icoideal}
(a) The set $I_H(A)$ is a two-sided ideal and a right $H$-coideal of~$T(X_H)$.

(b) The ideal $I_H(A)$ is preserved by all $H$-comodule algebra endomorphisms of~$T(X_H)$.
\end{Prop}

This result allows us to paraphrase the classical theory of rings with polynomial identities 
(see, e.g.,~\cite{Row}) by saying that $I_H(A)$ is a $T$-ideal.

\pf
(a) For each $\mu \in \RAlgComod{H}(T(X_H),A)$ the kernel $\Ker\mu$ is
a two-sided ideal of~$T(X_H)$. Hence so is the intersection~$I_H(A)$.

Let us show that $\Ker\mu$ is a right $H$-coideal of~$T(X_H)$ for each
comodule algebra morphism $\mu : T(X_H) \to A$. 
If we fix a basis $\{h_i\}_{i\in I}$ of~$H$, then there is a family of linear endomorphisms
$\{\delta_i\}_{i\in I}$ of~$T(X_H)$ such that 
the coaction of~$H$ on~$T(X_H)$ can be written
$$\delta(P) = \sum_{i\in I}\, \delta_i(P) \otimes h_i$$
for all~$P\in T(X_H)$.
If $P\in \Ker\mu$, then
$$\sum_{i\in I}\, \mu(\delta_i(P)) \otimes h_i
=  (\mu \otimes \id_H) \, \delta(P) = \delta (\mu(P)) = 0 \, .$$
It follows that $\mu(\delta_i(P)) = 0$ for all $i\in I$,
hence $\delta(P)$ belongs to $\Ker\mu \otimes H$.
We have thus established that $\Ker\mu$ is a right $H$-coideal of~$T(X_H)$.
Using a similar argument, one proves that 
the intersection $I_H(A)$ is a right $H$-coideal.

(b) If $f: T(X_H) \to T(X_H)$ is an $H$-comodule algebra morphism, then
so is  $\mu \circ f : T(X_H) \to A$ 
for all $\mu \in \RAlgComod{H}(T(X_H),A)$.
Hence, if $P\in I_H(A)$, then $\mu (f(P)) = 0$ for all $\mu \in \RAlgComod{H}(T(X_H),A)$,
which implies that $f(P)$ belongs to~$I_H(A)$.
\epf

The tensor algebra has a natural grading for which $T^r(X_H)$ ($r\geq 0$) 
is the subspace of all elements of~$T(X_H)$ of degree~$r$.
Equation~\eqref{TXcoaction2} implies that
\begin{equation}\label{Tgraded-comod}
\delta(T^r(X_H)) \subset T^r(X_H) \otimes H \, .
\end{equation}

\begin{Prop}\label{graded-ideal}
If $k$ is infinite, then the submodule $I_H(A)$ of~$T(X_H)$ is graded
with $I^r_H(A) = I_H(A) \bigcap T^r(X_H)$ for all $r \geq 0$.
\end{Prop}

\pf
Expand $P \in I_H(A)$ into a finite sum with homogeneous parts:
$$P = \sum_{r\geq 0}\, P_r \, ,$$
where $P_r \in T^r(X_H)$. It suffices to show that $P_r \in I_H(A)$ for all $r\geq 0$.

Given $\lambda \in k$, consider the algebra endomorphism $\lambda_*$ of~$T(X_H)$ 
defined by $\lambda(X_h) = \lambda X_h$ for all $h\in H$.
It is easy to check that $\lambda_*$ is an $H$-comodule map.
If $\mu : T(X_H) \to A$ is an $H$-comodule algebra morphism,
then so is $\mu \circ \lambda_*$. Consequently, 
$$\sum_{r\geq 0}\, \lambda^r \mu(P_r) = (\mu \circ \lambda_*)(P) = 0\, .$$
The $A$-valued polynomial $\sum_{r\geq 0}\, \lambda^r \mu(P_r)$ takes zero values for all $\lambda \in k$.
By the assumption on~$k$, this implies that its coefficients are all zero,
i.e., $\mu(P_r) = 0$ for all $r\geq 0$. 
Since this holds for all $\mu \in \RAlgComod{H}(T(X_H),A)$,
we obtain $P_r \in I_H(A)$ for all $r\geq 0$.
\epf

\begin{Prop}\label{degree0-identity}
Let $k$ be an infinite field.

(a) We have $I^0_H(A) = \{0\}$.

(b) If there is an injective comodule map $H \to A$, then $I^1_H(A)  = \{0\}$.
\end{Prop}

\pf
(a) The dimension of~$I^0_H(A)$ is at most one.
Since $\mu(1) = 1 \neq 0$ for all $\mu \in \RAlgComod{H}(T(X_H),A)$,
there cannot exist nonzero $H$-identities in degree zero.

(b) Any element of~$T(X_H)$ of degree one is of the form~$X_h$ for some $h\in H$.
Let $u: H\to A$ be an injective comodule map. 
It extends to a morphism $\mu \in \RAlgComod{H}(T(X_H),A)$ 
uniquely determined by $\mu(X_h) = u(h)$ for all $h\in H$.
If $X_h \in I_H(A)$, then $u(h) = \mu(X_h) = 0$.
By the injectivity of~$u$, we have $h = 0$, hence $X_h = 0$.
\epf

When the condition of Proposition~\ref{degree0-identity}\,(b) is not satisfied,
there may exist $H$-identities of degree one, 
as shown by Proposition~\ref{degree1-identity} below.

Let $H$ be a Hopf algebra and $A = k$ be the trivial one-dimensional $H$-comodule
algebra.
(In this case there are no injective comodule maps $H \to A$ unless $\dim H = 1$.)
Let us determine the ideal $I_H(k)$ of $H$-identities for~$k$.
It is well known (see~\cite[Sect.~3.3.1]{Abe}) that the vector space of 
$H$-comodule maps $H \to k$ coincides with the vector space of right integrals for~$H$
and that this vector space is either zero or one-dimensional.

If $H$ has no nonzero right integrals, then any element of the augmentation
ideal $T^+(X_H) = \bigoplus_{r\geq 1}\, T^r(X_H)$ is an $H$-identity for~$k$.
Hence,
\begin{equation}\label{T+}
I_H(k) = T^+(X_H)\, .
\end{equation}

Let us now consider the case where $H$ has a nonzero right integral.
This is the case for instance when $H$ is a finite-dimensional Hopf algebra.

\begin{Prop}\label{degree1-identity}
If $N: H \to k$ is a nonzero right integral, then
$$I^1_H(k) = \{ X_h \in T^1(X_H) \, |\, N(h) = 0 \} \, ,$$
and $I_H(A)$ is the two-sided ideal of $T(X_H)$ generated by~$I^1_H(k)$.
\end{Prop}

\pf
Let $\mu : T(X_H) \to k$ be an $H$-comodule algebra morphism. 
Its restriction to~$T^1(X_H) \cong H$ is an $H$-comodule map,
hence a right integral for~$H$, which
is necessarily a scalar multiple of~$N$. 
It follows that an element of~$T^1(X_H)$ is an $H$-identity for~$k$
if and only if it is in the kernel of~$N$. 

Let $\{h_i\}_{i\in I}$ be a basis of~$\Ker (N)$ and $h_0 \in H$ such that $N(h_0) = 1$.
Any element of $T(X_H)$ is a linear combination of monomials in the indeterminates
$X_{h_i}$ ($i\in I$) together with the indeterminate~$X = X_{h_0}$.
Clearly, $\mu$ vanishes on any monomial containing at least an indeterminate~$X_{h_i}$,
where $i \in I$. On the contrary, if $\lambda \neq 0$, then 
$\mu$ takes a nonzero value on any monomial that is a pure power of~$X$.
Therefore, $I_H(A)$ is the ideal generated by the indeterminates~$X_{h_i}$,
where $i\in I$.
\epf

\begin{Rems}\label{reduced tensor algebra}
(a) Instead of defining the $H$-identities in the tensor
algebra~$T(X_H)$, we could define them in the 
\emph{reduced tensor algebra}~$\bar{T}(X_H)$, 
which we define as the quotient of~$T(X_H)$
by the two-sided ideal generated by~$X_{1_H}  -1$, where $1_H$ is
the unit of the algebra~$H$. 
Certain identities would be simpler in~$\bar{T}(X_H)$, but we would 
lose the natural grading of~$T(X_H)$ 
as the ideal generated by~$X_{1_H}  -1$ is not graded.

(b) If we have a family $(A_{\alpha})_{\alpha}$ of $H$-comodule
algebras, we define an $H$-identity for the family
$(A_{\alpha})_{\alpha}$ to be an element $P\in T(X_H)$ that is an
$H$-identity for each $H$-comodule algebra~$A_{\alpha}$ of the family. 
We may then consider the intersection 
$\bigcap_{\alpha} I_H(A_{\alpha})$, 
which is a two-sided ideal and a right $H$-coideal of~$T(X_H)$. 
\end{Rems}

\subsection{Central polynomials}\label{central-identities}

Let $A$ be an $H$-comodule algebra. We say that $P \in T(X_H)$ is a
\emph{central polynomial for~$A$} if $PQ - QP$ is an $H$-identity for~$A$
for all $Q \in T(X_H)$. 
Since
$$ P(Q_1Q_2) - (Q_1Q_2P) = (PQ_1 - Q_1P)\, Q_2 + Q_1\, (PQ_2 - Q_2P)$$
for all $Q_1$, $Q_2 \in T(X_H)$, we conclude by induction on the degree
that $P$~is a central polynomial for~$A$ if and only if 
$$PX_h - X_hP \in I_H(A)$$
for all $h \in H$.
To check that $P$ is a central polynomial,
it is therefore enough to check that $\mu(P)$ commutes with $\mu(X_h)$ in~$A$ 
for all $\mu  \in \RAlgComod{H}(T(X_H),A)$ and all $h \in H$.

In \S~\ref{coinv-UH(A)} we shall establish under certain assumptions that
each element $P_{h^{[1]}, \ldots, h^{[r]}} \in T(X_H)$ of
Lemma~\ref{degree2-identity} is a central polynomial.

\section{The universal $H$-comodule algebra}\label{univ-comod-alg}

Let $A$ be an $H$-comodule algebra and 
$I_H(A)$ the two-sided ideal of $H$-identities for~$A$
defined by~\eqref{I(A)-def}.
Consider the quotient algebra
\begin{equation}\label{def-U}
\UU_H(A) = T(X_H)/I_H(A) \, .
\end{equation}

\begin{Prop}
(a) The algebra $\UU_H(A)$ has a unique structure of an $H$-comodule algebra
such that the projection $T(X_H) \to \UU_H(A)$ is an $H$-comodule algebra morphism.

(b) Any $\mu \in \RAlgComod{H}(T(X_H),A)$ induces an $H$-comodule algebra morphism 
$$\mu : \UU_H(A) \to A \, . $$
\end{Prop}

\pf
(a) This is an immediate consequence of Proposition~\ref{Icoideal}\,(a).

(b) This follows from the vanishing of~$\mu(I_H(A))$.
\epf

By definition of~$\UU_H(A)$, all $H$-identities for~$A$ are satisfied
(i.e., vanish) in~$\UU_H(A)$
and any quotient $H$-comodule algebra of~$T(X_H)$ in
which all $H$-identities for~$A$ are satisfied is a homomorphic image of~$\UU_H(A)$.
In this sense, $\UU_H(A)$ is the universal algebra in which all $H$-identities for~$A$
are satisfied.
We call $\UU_H(A)$ the \emph{universal $H$-comodule algebra} of~$A$.

By Proposition~\ref{graded-ideal}, if $k$ is infinite, then $I_H(A)$ is graded.
For each $r\geq 0$, let $\UU^r_H(A)$ be the image of $T^r(X_H)$ in~$\UU_H(A)$
under the projection map $T(X_H) \to \UU_H(A)$.
We have $1 \in \UU^0_H(A)$.

\begin{Prop}
If $k$ is infinite, then
the $H$-comodule algebra $\UU_H(A)$ is graded, namely
$$\UU_H(A) = \bigoplus_{r\geq 0}\, \UU^r_H(A) \, ,$$
$$\UU^r_H(A) \, \UU^s_H(A) \subset \UU^{r+s}_H(A) \quad \text{and} \quad
\delta(\UU^r_H(A)) \subset \UU^r_H(A) \otimes H$$
for all $r$, $s\geq 0$.
\end{Prop}

\pf
Clearly, $\UU_H(A) = \sum_{r\geq 0}\, \UU^r_H(A)$. To prove that this sum is direct,
it suffices to check that, if $\{u_r\}_{r\geq 0}$ is a family of elements of~$\UU_H(A)$
such that $u_r \in \UU^r_H(A)$ for all $r\geq 0$ and
$\sum_{r\geq 0}\, u_r = 0$, then $u_r = 0$ for all $r\geq 0$.
Represent each $u_r$ by an element $w_r \in T^r(X_H)$.
We have $\sum_{r\geq 0}\, w_r \in I_H(A)$. By Proposition~\ref{graded-ideal},
$I_H(A)$ is graded. This implies that $w_r \in I_H(A)$ for each~$r\geq 0$.
Therefore, $u_r = 0$ for each $r\geq 0$, as desired.

The inclusion $\UU^r_H(A) \, \UU^s_H(A) \subset \UU^{r+s}_H(A)$
follows from the fact that $T(X_H)$ is a graded algebra,
and the inclusion $\delta(\UU^r_H(A)) \subset \UU^r_H(A) \otimes H$
from~\eqref{Tgraded-comod}.
\epf

We now use Proposition~\ref{degree1-identity} to give a simple example
for which we can compute~$\UU_H(A)$. In this example, 
$A$ is the trivial $H$-comodule algebra~$k$.

\begin{Prop}
(a) If there are no nonzero right integrals for~$H$, then $\UU_H(k) \cong k$.

(b) If there is a nonzero right integral for~$H$, then $\UU_H(k)$
is isomorphic as an $H$-comodule algebra to~$k[X]$ 
with trivial coaction.
\end{Prop}

\pf
(a) This is a consequence of~\eqref{T+}.

(b) The isomorphism $\UU_H(k) \cong k[X]$ follows from the proof of
Proposition~\ref{degree1-identity}.
It remains to show that the coaction on~$k[X]$ is given by $\delta(X) = X \otimes 1$. 
We resume the notation of the second part of the proof of
Proposition~\ref{degree1-identity}.
Expand
$$\Delta(h_0) = h_0 \otimes x_0 + \sum_{i\in I}\, h_i \otimes x_i \, ,$$
where $x_0$ and $\{x_i\}_{i \in I}$ are elements of~$H$.
Since $N : H \to k$ is a comodule map, 
$$N(h_0) \otimes x_0 + \sum_{i\in I}\, N(h_i) \otimes x_i 
= (N \otimes \id_H) \bigl(\Delta(h_0)\bigr)
= N(h_0)\, 1 \otimes 1\, .$$
Replacing $N(h_i)$ by its value, we obtain
$1 \otimes x_0 = 1 \otimes 1$, hence $x_0 = 1$.
Therefore, by definition of the coaction of~$H$ on~$T(X_H)$,
$$\delta(X) = X \otimes 1 + \sum_{i\in I}\, X_{h_i} \otimes x_i 
\in X \otimes 1 + I_H(k) \otimes H \, .$$
We obtain the desired conclusion by passing to the quotient by~$I_H(k)$.
\epf

We end this section with a few questions on the structure of~$\UU_H(A)$. 

(i) If $A$ is free as a module over the subalgebra of coinvariants,
does the same hold for~$\UU_H(A)$?

(ii) If $A$ is free as a module over its center,
is $\UU_H(A)$ free over its center~$\ZZ_H(A)$,
possibly after some localization of the latter?

(iii) Is a suitable central localization of~$\UU_H(A)$ a flat deformation of~$A$
over the center~$\ZZ_H(A)$? 
How is $Z(A) \otimes_{\ZZ_H(A)} \UU_H(A)$ related to~$A$, 
where $Z(A)$ is the center of~$A$?

In the sequel we shall provide answers to these questions 
under some hypotheses.

\section{Detecting $H$-identities}\label{universal evaluation}

We again fix a Hopf algebra~$H$.
The aim of this section is construct 
an $H$-comodule algebra morphism whose domain is~$T(X_H)$ and 
whose kernel is~$I_H(A)$ for any $H$-comodule algebra~$A$
satisfying the following assumptions.

\subsection{The assumptions}\label{assumptions}

Let $A$ be an $H$-comodule algebra and 
$B = A^H$ be the subalgebra of coinvariant elements of~$A$.
We make the following assumptions on~$A$:

(i) $B$ is central in~$A$ and

(ii) there is a $B$-linear $H$-comodule isomorphism $u: B\otimes H \to A$, 
where the coaction on~$B\otimes H$ is equal to~$\id_B \otimes \Delta$. 

In Section~\ref{universal cocycle} we shall give examples of $H$-comodule algebras
for which these assumptions hold.

\subsection{The symmetric algebra}\label{S(t_H)}

Let $t_H$ be another copy of the underlying vector space of~$H$. 
We denote the identity map from $H$ to~$t_H$ by
$h\mapsto t_h$ ($h\in H$).

Let $S(t_H)$ be the {\em symmetric algebra} over the $k$-vector space~$t_H$.
It is a quotient of the tensor algebra~$T(t_H)$ and inherits the grading of the latter.
For each $r\geq 0$, let us denote $S^r(t_H)$ the subspace of elements of degree~$r$.
We have
$$S(t_H) = \bigoplus_{r\geq 0}\, S^r(t_H)$$
with $S^0(t_H) = k$, and $S^1(t_H) = t_H$.
If $\{h_i\}_{i\in I}$ is a basis of~$H$, then $S(t_H)$ is
isomorphic to the polynomial algebra over the indeterminates $\{t_{h_i}\}_{i\in I}$.

The linear isomorphism $t : H \to t_H$ extends to an injective linear
map $H \to S(t_H)$, which we still denote~$t$ and whose image is~$S^1(t_H)$.
The algebra $S(t_H)$ satisfies the following well-known universal property:
for any commutative algebra~$R$ and any linear map $g : H \to R$,
there is a unique algebra morphism $f : S(t_H) \to R$ such that $f\circ t = g$.

\subsection{The universal evaluation morphism}\label{univ-eval-map}

Given an algebra~$A$,
we endow the tensor product $S(t_H) \otimes A$ with its natural algebra structure.
The grading of~$S(t_H)$ induces a grading on the algebra $S(t_H) \otimes A$:
for each $r\geq 0$,
the subspace of elements of degree~$r$ in $S(t_H) \otimes A$ is equal to~$S^r(t_H) \otimes A$.

If, in addition, $A$ is an $H$-comodule algebra, then
$S(t_H) \otimes A$ is an $H$-comodule algebra with the $S(t_H)$-linear coaction
$$\id_{S(t_H)} \otimes \delta : S(t_H) \otimes A \to S(t_H) \otimes A \otimes H\, ,$$
where $\delta : A \to A \otimes H$ is the coaction of~$A$.

Suppose now that $A$ satisfies the assumptions of \S~\ref{assumptions}.
Using the $H$-comodule algebra~$T(X_H)$ of Section~\ref{identities},
we consider the algebra morphism
$$\mu_{A} : T(X_H) \to S(t_H) \otimes A$$
defined for all $h\in H$ by
\begin{equation}\label{universal map}
\mu_{A}(X_h) = \sum_{(h)}\, t_{h\sw1} \otimes u(1_B \otimes h\sw2) \, ,
\end{equation}
where $u: B\otimes H \to A$ is the isomorphism of~\S~\ref{assumptions}.

\begin{Lem}\label{lem-mu-graded}
(a) The algebra morphism $\mu_{A} : T(X_H) \to S(t_H) \otimes A$ 
sends $T^r(X_H)$ into $S^r(t_H) \otimes A$ for each $r\geq 0$.

(b) The morphism $\mu_{A}$ is injective on 
$T^0(X_H) \oplus T^1(X_H) = k \oplus X_H$ 
and its kernel is contained in $\bigoplus_{r\geq 2}\, T^r(X_H)$.
\end{Lem}

\pf
(a) This follows from \eqref{universal map}.

(b) Since $\mu_{A}(1) = 1 \otimes u(1_B \otimes 1_H) \neq 0$,
the morphism $\mu_{A}$ is injective on~$T^0(X_H)$. 
It is injective on~$T^1(X_H) = X_H$ since
$$\bigl( (\id_{S(t_H)}\otimes \eps)\circ (\id_B \otimes u^{-1}) \circ \mu_{A}\bigr)(X_h)
= \sum_{(h)}\, t_{h\sw1} \otimes \eps(h\sw2) 1_B = t_h \otimes 1_B$$
for all $h\in H$.
\epf

The morphism $\mu_{A}$ possesses the following important properties.

\begin{Lem}\label{lem-mu-univ}
(a) The morphism $\mu_{A} : T(X_H) \to S(t_H) \otimes A$ is an
$H$-comodule algebra morphism.

(b) For every $\mu \in \RAlgComod{H}(T(X_H),A)$ there is a unique algebra morphism 
$$\lambda : S(t_H) \to B$$ 
such that $\mu = (\lambda \otimes \id_A) \circ \mu_{A}$.
\end{Lem}

In view of the previous lemma, we call $\mu_{A}$ the {\em universal
evaluation morphism} for~$A$.

\pf 
(a) It is enough to check that the restriction of $\mu_A$ to~$X_H$
is an $H$-comodule map. Now, $\mu_A = (\id_{S(t_H)} \otimes u) \circ
\mu'_A$, where 
$$\mu'_A : X_H \to S(t_H) \otimes B \otimes H$$
is the linear map defined by 
$\mu'_A(X_h) = \sum_{(h)}\, t_{h\sw1} \otimes 1_B \otimes h\sw2$ 
for all $h \in H$. 
Using the definition of the coaction on $X_H$ and on $S(t_H) \otimes B \otimes H$,
and the coassociativity of the comultiplication of~$H$, one checks easily that
$\mu'_A$ is a comodule map. 
Since $u$ also is a comodule map, the conclusion follows.

(b) Since we are dealing with algebra morphisms, it suffices to consider
their restriction to~$X_H$. Moreover, $u : B\otimes H \to A$ being
$B$-linear, it is enough to prove that for every $\mu \in
\RComod{H}(X_H,B\otimes H)$ there is a unique linear map $\lambda :
t_H \to B$ such that $\mu = (\lambda \otimes \id_B\otimes \id_H)
\circ \mu'_{A}$.

We first prove the uniqueness of~$\lambda$. 
If $\mu = (\lambda \otimes \id_B\otimes \id_H) \circ \mu'_{A}$, then 
$(\id_B \otimes \eps) \circ \mu = \lambda \otimes \id_B\otimes \eps$. 
Hence, for all $h\in H$,
\begin{eqnarray*}
(\id_B \otimes \eps)\bigl(\mu (X_h)\bigr)
& = & \sum_{(h)}\, \lambda(t_{h\sw1})  1_B \, \eps(h\sw2) \\
& = & \sum_{(h)}\, \lambda(t_{h\sw1\eps(h\sw2)})  
=  \lambda(t_h) \, .
\end{eqnarray*}
This shows that $\lambda$ is uniquely determined.

Let us next prove the existence of~$\lambda$. 
For $\mu \in \RAlgComod{H}(T(X_H),B\otimes H)$ and $h\in H$, 
expand~$\mu(X_h)$ in~$B\otimes H$ in the form
$$\mu(X_h) = \sum_{j}\, \mu_{j}(h) \otimes \nu_{j}(h) \, ,$$
where $\mu_{j}(h) \in B$ and $\nu_{j}(h) \in H$.
If we set $\lambda(t_h) = \sum_{j}\, \mu_{j}(h) \, \eps(\nu_{j}(h)) \in B$,
then
\begin{eqnarray*}
\bigl((\lambda \otimes \id_B\otimes \id_H) \circ \mu'_{A}\bigr)
(X_h) & = & \sum_{(h)}\, \lambda(t_{h\sw1}) \otimes h\sw2 \\
& = & \sum_{j, \, (h)}\, \mu_{j}({h_{(1)}}) \, \eps(\nu_{j}({h_{(1)}})) \otimes h_{(2)} \\
& = & \sum_{j, \, (\nu_j(h))}\, \mu_{j}(h) \otimes \eps(\nu_{j}(h)_{(1)}) \, \nu_{j}(h)_{(2)} \\
& = & \sum_{j}\, \mu_{j}(h) \otimes \nu_{j}(h) 
= \mu(X_h)
\end{eqnarray*}
for all $h\in H$.
The third and fourth equalities follow respectively from the facts 
that $\mu$ is a comodule map
and $\eps$ is the counit. 
\epf

Part~(b) of the previous proof can be summarized by the following
sequence of natural bijections:
\begin{eqnarray*}
\RAlgComod{H}(T(X_H),A)
& \cong & \RComod{H}(H,A) \\
& \cong & \RComod{H}(H,B\otimes H)\\
& \cong & \Hom(H,B) \\
& \cong & \Alg(S(t_H), B) \, ,
\end{eqnarray*}
where the last bijection holds because the algebra~$B$ is commutative.

We now state the main result of this section: it implies that the
$H$-identities for~$A$ are detected by the morphism~$\mu_A$.

\begin{Thm}\label{thm-mu-kernel}
Let $A$ be an $H$-comodule algebra satisfying 
the assumptions of \S~\ref{assumptions}.

(a) We have $\Ker \bigl(\mu_A : T(X_H) \to S(t_H)\otimes A \bigr) \subset I_H(A)$.

(b) If the ground field $k$ is infinite, then
$$I_H(A) = \Ker \bigl(\mu_A : T(X_H) \to S(t_H)\otimes A \bigr) \, .$$
\end{Thm}

\pf 
(a) Let $P\in T(X_H)$ such that $\mu_A(P) = 0$. 
By Lemma~\ref{lem-mu-univ}\,(b), 
any morphism $\mu \in \RAlgComod{H}(T(X_H),A)$ is
of the form $\mu = (\lambda \otimes \id_A) \circ \mu_A$ for some
$\lambda \in \Alg(S(t_H), B)$. Therefore, $\mu(P) = (\lambda \otimes
\id_H) \bigl(\mu_A(P)\bigr) = 0$. Hence, $P \in I_H(A)$.

(b) By Part~(a) it is enough to check that any $P \in I_H(A)$ belongs to the
kernel of~$\mu_A$. 
Since $\id_{S(t_H)} \otimes u$ is a linear isomorphism, it suffices to check
that $\mu'_A(P) = 0$ for any $P\in I_H(A)$, where 
$\mu'_A : X_H \to S(t_H) \otimes B \otimes H$ is the map introduced in the proof
of Lemma~\ref{lem-mu-univ}.
Expanding $\mu'_A(P)$ in a basis $\{h_i\}_{i\in I}$ of~$H$, we have
$$\mu'_A(P) = \sum_{i\in I}\, \mu^{(i)}_A(P) \otimes h_i  \in S(t_H) \otimes B \otimes H \, ,$$
where $\mu^{(i)}_A(P) \in S(t_H) \otimes B$ for all $i\in I$. 
By Lemma~\ref{lem-mu-univ}\,(b),
if $P \in I_H(A)$, then 
$$\sum_{i\in I}\, \lambda \bigl(\mu^{(i)}_A(P) \bigr) \otimes h_i 
= (\lambda \otimes \id_{B\otimes H}) \bigl(\mu'_A(P)\bigr) = 0$$ 
for all $\lambda \in \Alg(S(t_H), B)$.
Therefore, for each $i\in I$, we have $\lambda(\mu^{(i)}_A(P)) \nolinebreak = \nolinebreak 0$
for all algebra morphisms $\lambda : S(t_H) \to B$. This means that the
$B$-valued polynomial $\mu^{(i)}_A(P)$ takes only zero values.
Since $k$ is an infinite field, this implies that $\mu^{(i)}_A(P) = 0$ for all $i\in I$.
Hence, $\mu'_A(P) = 0$. 
\epf

We have the following consequence of Lemma~\ref{lem-mu-graded} and of
Theorem~\ref{thm-mu-kernel} 
for the universal $H$-comodule algebra~$\UU_H(A)$ defined by~\eqref{def-U}.

\begin{Cor}\label{inject}
Let $k$ be an infinite field and 
$A$ be an $H$-comodule algebra satisfying the assumptions of
\S~\ref{assumptions}.

(a) The universal evaluation morphism $\mu_A$ induces an injective
$H$-co\-mod\-ule algebra morphism
$$\mu_A :
\UU_H(A) = T(X_H)/I_H(A) \hookrightarrow S(t_H) \otimes A \, .$$

(b) We have
$\UU_H^0(A) \cong k$ and $\UU_H^1(A) \cong X_H$.
\end{Cor}

The previous corollary has an interesting consequence,
which guarantees the existence of nonzero $H$-identities for~$A$
in the finite-dimensional case.

\begin{Cor}\label{ident-exist}
Let $k$ be an infinite field and 
$A$ be an $H$-comodule algebra satisfying the assumptions of
\S~\ref{assumptions}.
If $H$ and $A$ are finite-dimensional and $H \neq k$, then $I_H(A) \neq 0$.
\end{Cor}

\pf
If $I_H(A) = 0$, then Corollary~\ref{inject} implies that
$T^r(X_H)$ injects into $S^r(t_H) \otimes A$. Set $n = \dim H$ and $d = \dim A$.
We then have
$$n^r = \dim T^r(X_H) \leq \dim (S^r(t_H) \otimes A) = d \, \binom{r+n-1}{n-1}\, ,$$
which is impossible for $n>1$ and large~ $r$.
\epf

\begin{Rem}
Suppose that the map $u: B\otimes H \to A$ of \S~\ref{assumptions}
preserves the units.
Since the universal evaluation morphism $\mu_{\alpha}$ sends $X_{1_H}$ to
the central element $t_1 \otimes u(1_B\otimes 1_H) = t_1 \otimes 1_A$, 
it induces an $H$-comodule algebra morphism $\bar{\mu}_{\alpha}: \bar{T}(X_H)
\to \bar{S}(t_H) \otimes A$, where $\bar{T}(X_H)$ is the
reduced tensor algebra of Remark~\ref{reduced tensor algebra}\,(a)
and $\bar{S}(t_H) = S(t_H)/(t_{1_H} - 1)$. As can be seen from the
computations in Section~\ref{Sweedler}, using the morphism $\bar{\mu}_{\alpha}$ 
would simplify certain formulas, but these formulas would
no longer be homogeneous with respect to the gradings 
of $T(X_H)$ and of~$S(t_H)$.
\end{Rem}

\subsection{The coinvariants of~$\UU_H(A)$}\label{coinv-UH(A)}

The following proposition shows how we can use~$\mu_A$ to detect 
the co\-\"\i n\-var\-i\-ant elements of~$\UU_H(A)$.

\begin{Prop}\label{coinv-UU}
Let $k$ be an infinite field and 
$A$ be an $H$-comodule algebra satisfying the assumptions of
\S~\ref{assumptions}.
Then 
$$\UU_H(A)^H = \{P \in \UU_H(A)\, |\, \mu_A(P) \in S(t_H) \otimes B \} \, .$$
\end{Prop}

\pf
This follows from Lemma~\ref{inject-coinv}, Corollary~\ref{inject}\,(a), 
and from the easily checked fact that
the subalgebra of coinvariants of $S(t_H) \otimes A$ is~$S(t_H) \otimes B$.
\epf

We give an interesting consequence of Proposition~\ref{coinv-UU}
for the elements $P_{h^{[1]}, \ldots, h^{[r]}} \in T(X_H)$ introduced in
Lemma~\ref{degree2-identity}.

\begin{Cor}\label{ex-central-identity}
Under the hypotheses of Proposition~\ref{coinv-UU},
$P_{h^{[1]}, \ldots, h^{[r]}}$
is a central polynomial for~$A$ for all $h^{[1]}, \ldots, h^{[r]} \in H$.
\end{Cor}

\pf
By Lemma~\ref{degree2-identity}, $P_{h^{[1]}, \ldots, h^{[r]}}$ is coinvariant.
Thus, by Proposition~\ref{coinv-UU}, 
its image $\mu_A(P_{h^{[1]}, \ldots, h^{[r]}})$ in~$S(t_H) \otimes A$ 
belongs to~$S(t_H) \otimes B$. 
Since $S(t_H) \otimes B$ is central in~$S(t_H) \otimes A$, the element
$\mu_A(P_{h^{[1]}, \ldots, h^{[r]}})$ commutes with $\mu_A(X_h)$ for all~$h\in H$. 
By Lemma~\ref{lem-mu-univ}\,(b),
$\mu(P_{h^{[1]}, \ldots, h^{[r]}})$ then commutes with~$\mu(X_h)$ 
for all $\mu \in \RAlgComod{H}(T(X_H),A)$ and all~$h\in H$. 
This shows that $P_{h^{[1]}, \ldots, h^{[r]}}$ is a central polynomial.
\epf

\section{The universal cocycle attached to a twisted algebra}\label{universal cocycle}

Let $H$ be a Hopf algebra.

\subsection{Twisted algebras}\label{Twisted products}

Consider the following general construction.
Let $\alpha : H \times H \to B$ be a convolution invertible bilinear map
with values in a commutative algebra~$B$.
Consider the \emph{twisted algebra} $B\otimes {}^{\alpha} H$,
which is the vector space~$B\otimes H$ equipped with the product~$\cdot_{\alpha}$ given by
\begin{equation}\label{twisted-multiplication}
(b\otimes x) \cdot_{\alpha} (c\otimes y)
= \sum_{(x), (y)}\, bc \, \alpha(x\sw1, y\sw1) \otimes  x\sw2 y\sw2
\end{equation}
for all $b$, $c\in B$ and $x$, $y \in H$.
If $\alpha$ is a two-cocycle of~$H$, 
which means that it satisfies the equations
\begin{equation}\label{cocycle}
\sum_{(x), (y)}\alpha(x\sw1,y\sw1)\, \alpha(x\sw2y\sw2,z)
= \sum_{(y), (z)} \alpha(y\sw1,z\sw1)\, \alpha(x,y\sw2z\sw2)
\end{equation}
for all $x,y,z \in H$,
then the product $\cdot_{\alpha}$ is associative.
If in addition $\alpha$ is normalized, i.e., if
\begin{equation}\label{normalized}
\alpha(x,1)  = \alpha(1,x) = \varepsilon(x)\, 1
\end{equation}
for all $x\in H$,
then $1\otimes 1$ is a unit for the product~$\cdot_{\alpha}$.

The algebra $A = B\otimes {}^{\alpha} H$ is an $H$-comodule algebra
with coaction 
$$\delta = \id_B \otimes \Delta : A = B\otimes H \to A \otimes H = B\otimes H \otimes H \, .$$ 
Its subalgebra~$A^H$ of coinvariants coincides with~$B\otimes k 1$. 
It is easy to check from~\eqref{twisted-multiplication}
that $A^H$ lies in the center of~$A$.
The assumptions of \S~\ref{assumptions} are satisfied for~$A$
with $B=A^H$ and~$u = \id_{B\otimes H}$.

The twisted algebras of the form~$B\otimes {}^{\alpha} H$ considered here coincide
with the central \emph{cleft $H$-Galois extensions} of~$B$,
see \cite{BCM}, \cite{DT1}, \cite[Prop.~7.2.3]{M2}.
(By central we mean that $B$ lies in the center of the extension.)

\subsection{The two-cocycle $\sigma$}\label{universal cocycle1}

Let $\alpha : H \otimes H \to k$ be
a normalized convolution-invertible two-cocycle
with values in the ground field~$k$
and $k\otimes {}^{\alpha} H$ be the corresponding twisted algebra.
To simplify notation, we henceforth denote $k\otimes {}^{\alpha} H$
by~${}^{\alpha} H$.

To the two-cocycle~$\alpha$ we associate a two-cocycle
$$\sigma : H \times H \to K_H$$
taking values in the field of fractions~$K_H$ of~$S(t_H)$
and defined by
\begin{equation}\label{sigma-def}
\sigma(x,y) = \sum_{(x), (y)}\, t_{x_{(1)}} t_{y_{(1)}} \,
\alpha(x_{(2)},y_{(2)}) \, t^{-1}_{x_{(3)} y_{(3)}}
\end{equation}
for all $x,y \in H$. 
In \eqref{sigma-def} we have used the 
the map $t^{-1} : H \to K_H$ uniquely defined
by the equations~\eqref{tbareqn} of Appendix~\ref{tbar-proof}.

By~\cite[Chap.~7]{M2}, $\sigma$ is a
two-cocycle on~$H$ and it is cohomologous to~$\alpha$ in~$K_H$. 
We call $\sigma$ the \emph{universal cocycle attached to}~$\alpha$; 
it can be viewed as the most general two-cocycle on~$H$ 
that is cohomologous to~$\alpha$.

Since $\alpha$ is convolution invertible, so is~$\sigma$. 
Its convolution inverse $\sigma^{-1}$ is given for all $x,y \in H$ by
\begin{equation}\label{sigma^{-1}-def}
\sigma^{-1}(x,y) = \sum_{(x), (y)}\, t_{x_{(1)} y_{(1)}} \, 
\alpha^{-1}(x_{(2)},y_{(2)}) \, t^{-1}_{x_{(3)}} t^{-1}_{y_{(3)}}\, ,
\end{equation}
where $\alpha^{-1}$ is the convolution inverse of~$\alpha$.

The two-cocycle $\sigma$ is almost normalized in the following sense.

\begin{Lem}\label{sigma(1,h)}
For all $x\in H$,
$$\sigma(1,x) = \sigma(x,1) = \eps(x) \, t_1
\quad \text{and} \quad \sigma^{-1}(1,x) = \sigma^{-1}(x,1) = \eps(x) \, t_1^{-1} \, .$$
\end{Lem}

\pf By~\eqref{normalized},~\eqref{sigma-def}, and~\eqref{tbareqn},
\begin{eqnarray*}
\sigma(1,x)
& = & \sum_{(x)}\, t_1 t_{x_{(1)}} \, \alpha(1,x_{(2)}) \, t^{-1}_{x_{(3)}} \\
& = & \sum_{(x)}\, t_1 t_{x_{(1)}} \, \eps(x_{(2)}) \, t^{-1}_{x_{(3)}} \\
& = & \sum_{(x)}\, t_1 t_{x_{(1)}} \, t^{-1}_{x_{(2)}} = \eps(x) \, t_1 \, .
\end{eqnarray*}
The other computations are similar.
\epf

\subsection{Two algebras of rational fractions}\label{universal cocycle2}

Let $\alpha : H \otimes H \to k$ and $\sigma : H \otimes H \to K_H$ be as
in \S~\ref{universal cocycle1}.

We define $S(t_H)_{\Theta}$ to be the subalgebra of~$K_H$ generated
by $t_H = t(H)$ and $t^{-1}_H = t^{-1}(H)$,
where $t^{-1} : H \to K_H$ is the map defined in Appendix~\ref{tbar-proof}.

This algebra depends only on the coalgebra structure of~$H$.
A description of~$S(t_H)_{\Theta}$ as an explicit localization of~$S(t_H)$
is given in Appendix~\ref{free-comm-coalg}.

By definition of the universal two-cocycle~$\sigma$, both $\sigma$ and 
its convolution inverse~$\sigma^{-1}$ take values in~$S(t_H)_{\Theta}$.

We recover the cocycle~$\alpha$ from the universal cocycle~$\sigma$ 
as follows.

\begin{Lem}\label{augmentation}
There is a unique algebra morphism $\epsilon : S(t_H)_{\Theta} \to k$ such that
$$\epsilon(t_x) = \epsilon(t^{-1}_x) = \varepsilon(x)$$
for all $x\in H$. 
Moreover, $\epsilon\bigl( \sigma(x,y) \bigr) = \alpha(x,y)$
and $\epsilon\bigl( \sigma^{-1}(x,y) \bigr) = \alpha^{-1}(x,y)$
for all $x$, $y\in H$.
\end{Lem}

\pf 
There is a unique algebra morphism $\epsilon : S(t_H) \to k$ such that
$\epsilon(t_x) = \varepsilon(x)$ for all $x\in H$. In order to check that
this morphism extends to~$S(t_H)_{\Theta}$, it suffices by~\eqref{tbareqn}
to check that
\begin{equation}\label{tbareqn01}
\sum_{(x)}\, \epsilon(t_{x\sw1}) \, \epsilon(t^{-1}_{x\sw2}) 
= \sum_{(x)}\, \epsilon(t^{-1}_{x\sw1})\, \epsilon(t_{x\sw2}) = \eps(x) \, 1 \, .
\end{equation}
for all $x\in H$.
Now,
\begin{equation*}
\sum_{(x)}\, \epsilon(t_{x\sw1}) \, \epsilon(t^{-1}_{x\sw2}) 
= \sum_{(x)}\, \varepsilon({x\sw1}) \, \varepsilon({x\sw2})
= \varepsilon\Bigl( \sum_{(x)}\, {x\sw1} \, {x\sw2} \Bigr)
= \eps(x) \, 1 \, .
\end{equation*}
The other equality is proved in a similar way.

We also have
\begin{eqnarray*}
\epsilon\bigl( \sigma(x,y) \bigr) 
& = & \epsilon\Bigl( \sum_{(x), (y)}\, t_{x_{(1)}} t_{y_{(1)}} \,
\alpha(x_{(2)},y_{(2)}) \, t^{-1}_{x_{(3)} h_{(y)}} \Bigr) \\
& = & \sum_{(x), (y)}\, \epsilon(t_{x_{(1)}}) \, \epsilon(t_{y_{(1)}}) \,
\alpha(x_{(2)},y_{(2)}) \, \epsilon(t^{-1}_{x_{(3)} y_{(3)}}) \\
& = & \sum_{(x), (y)}\, \eps({x_{(1)}}) \, \eps({y_{(1)}}) \,
\alpha(x_{(2)},y_{(2)}) \, \eps({x_{(3)} y_{(3)}}) \\
& = & \sum_{(x), (y)}\, \eps({x_{(1)}}) \, \eps({y_{(1)}}) \,
\alpha(x_{(2)},y_{(2)}) \, \eps({x_{(3)}) \, \eps(y_{(3)}}) \\
& = & \alpha(x,y) \, .
\end{eqnarray*}
The computation of~$\epsilon( \sigma^{-1}(x,y))$ is done in a similar way.
\epf

There is a smaller algebra in which both $\sigma$ and~$\sigma^{-1}$ take values,
namely the subalgebra $\BB_H^{\alpha}$ of~$K_H$
generated by the values of $\sigma$ and of~$\sigma^{-1}$.
The algebra $\BB_H^{\alpha}$ depends on the two-cocycle~$\alpha$ and 
is a subalgebra of~$S(t_H)_{\Theta}$:
$$\BB_H^{\alpha} \subset S(t_H)_{\Theta} \subset K_H\, .$$

\begin{Prop}\label{Bnozerodivisors}
The algebras $\BB_H^{\alpha}$ and $S(t_H)_{\Theta}$ have the following properties.

(a) They are augmented.

(b) They have no zero-divisors.

(c) If $H$ is finite-dimensional, then they are finitely generated.
\end{Prop}

\pf (a) This follows from Lemma~\ref{augmentation}.

(b) This holds since both algebras lie in the field~$K_H$.

(c) If $\{h_i\}_i$ is a basis of~$H$, then
the algebra $\BB_H^{\alpha}$ (resp.\ $S(t_H)_{\Theta}$)
is generated by the elements $\sigma(h_i,h_j)$ and~$\sigma^{-1}(h_i, h_j)$
(resp.\ by the elements $t_{h_i}$ and~$t^{-1}_{h_i}$)
If $H$ is finite-dimensional, these generators are in finite number. 
\epf

It is natural to ask how $\BB_H^{\alpha}$ depends on the two-cocycle~$\alpha$.
In particular, if two arbitrary two-cocycles $\alpha, \beta : H\times H \to k$ are cohomologous,
do we have $\BB_H^{\alpha} = \BB_H^{\beta}$ inside~$S(t_H)_{\Theta}$?
We do not know the answer to this question in general.
Nevertheless, there is a positive answer in the following case.

\begin{Prop}\label{dependanceBB}
Suppose that $\alpha, \beta : H\times H \to k$ are normalized con\-vo\-lu\-tion-invertible
two-cocycles such that there is a convolution-invertible linear form $\lambda : H \to k$ 
verifying the following two conditions for all $x$, $y\in H$:
\begin{equation}\label{beta-cond1}
\beta(x,y) = \sum_{(x), (y)}\, \lambda({x_{(1)}}) \lambda({y_{(1)}}) \,
\alpha(x_{(2)},y_{(2)}) \, \lambda^{-1}({x_{(3)}) y_{(3)}})\, ,
\end{equation}
where $\lambda^{-1}$ is the convolution inverse of~$\lambda$, and
\begin{equation}\label{beta-cond2}
\sum_{(x)}\, \lambda(x_{(1)}) x_{(2)}  =  \sum_{(x)}\,  x_{(1)} \lambda(x_{(2)}) \, .
\end{equation}
Then $\BB_H^{\alpha} = \BB_H^{\beta}$.
\end{Prop}

Condition~\eqref{beta-cond1} states that $\alpha$ and $\beta$ are 
cohomologous two-cocyles.
Condition~\eqref{beta-cond2} states that the linear form~$\lambda$ 
is \emph{lazy} in the sense of~\cite[\S~1]{BC}; 
this condition is automatically satisfied when $H$ is cocommutative.

\pf
Let $\sigma$ and $\sigma' : H\times H \to S(t_H)_{\Theta}$ be the 
universal cocycles attached to $\alpha$ and~$\beta$, respectively.
By definition (see~\eqref{sigma-def}),
$$\sigma'(x,y) = \sum_{(x), (y)}\, t_{x_{(1)}} t_{y_{(1)}} \,
\beta(x_{(2)},y_{(2)}) \, t^{-1}_{x_{(3)} y_{(3)}}\, .$$
Let us express $\sigma'(x,y)$ in terms of the elements~$\sigma(x,y)$
given by~\eqref{sigma-def}.
For all $x$, $y\in H$, we have
\begin{eqnarray*}
\sigma'(x,y)
& = & \sum_{(x), (y)}\, t_{x_{(1)}} t_{y_{(1)}} \,
\beta(x_{(2)},y_{(2)}) \, t^{-1}_{x_{(3)} y_{(3)}} \\
& = & \sum_{(x), (y)}\, t_{x_{(1)}} t_{y_{(1)}} \,
\lambda({x_{(2)}}) \lambda({y_{(3)}}) \,
\alpha(x_{(3)},y_{(3)}) \, \lambda^{-1}({x_{(4)}) y_{(4)}}) \, t^{-1}_{x_{(5)} y_{(5)}} \\
& = & \sum_{(x), (y)}\, \lambda({x_{(1)}}) \lambda({y_{(1)}}) \,
t_{x_{(2)}} t_{y_{(2)}} \,
\alpha(x_{(3)},y_{(3)})  \, t^{-1}_{x_{(4)} y_{(4)}} \, \lambda^{-1}({x_{(5)}) y_{(5)}}) \\
& = & \sum_{(x), (y)}\, \lambda({x_{(1)}}) \lambda({y_{(1)}}) \,
\sigma(x_{(2)},y_{(2)})  \, \lambda^{-1}({x_{(3)}) y_{(3)}}) \, .
\end{eqnarray*}
The third equality follows from Condition~\eqref{beta-cond2} on~$\lambda$ 
and from the same condition
on~$\lambda^{-1}$ (the latter is an immediate consequence of the former).
Denoting $\sigma'{}^{-1}$ the convolution inverse of~$\sigma'$, we similarly obtain
\begin{equation*}
\sigma'{}^{-1}(x,y)
= \sum_{(x), (y)}\, \lambda^{-1}({x_{(1)}}) \lambda^{-1}({y_{(1)}}) \,
\sigma^{-1}(x_{(2)},y_{(2)})  \, \lambda({x_{(3)}) y_{(3)}})
\end{equation*}
for all $x$, $y\in H$.
Thus, the generators $\sigma'{}(x,y)$ and $\sigma'{}^{-1}(x,y)$ of~$\BB_H^{\beta}$
are linear combinations of the generators of~$\BB_H^{\alpha}$.
It follows that
$\BB_H^{\beta} \subset \BB_H^{\alpha}$.
Exchanging the roles of $\alpha$ and~$\beta$, we obtain the reverse inclusion.
\epf

\begin{Rems}
(i) The transcendence degree of the field of fractions of~$\BB_H^{\alpha}$ over~$k$
cannot exceed the dimension of~$H$ over~$k$.

(ii) In the case where $H$ is finite-dimensional,
one can obtain a presentation of~$\BB_H^{\alpha}$ by generators and relations 
using standard monomial order techniques of commutative algebra 
(see, e.g.,~\cite[Chap.~15]{Eis}). 
In Section~\ref{Sweedler}
we shall give a presentation of~$\BB_H^{\alpha}$ for Sweedler's Hopf algebra.
\end{Rems}

\subsection{The integrality of $S(t_H)_{\Theta}$ over $\BB_H^{\alpha}$}\label{integrality}

We end this section by asking the following question:
under which conditions on the Hopf algebra~$H$ and on the
two-cocycle~$\alpha$ is $S(t_H)_{\Theta}$ integral over~$\BB_H^{\alpha}$?
If the answer to this question is positive, 
then the field of fractions~$K_H$ of~$S(t_H)_{\Theta}$ is an algebraic extension 
of the field of fractions~$\KK_H^{\alpha}$ of~$\BB_H^{\alpha}$.
The integrality of~$S(t_H)_{\Theta}$ over~$\BB_H^{\alpha}$
will be needed in Theorem~\ref{formthm} below.

We do not know how to answer this question in general.
We do not even have an example where $S(t_H)_{\Theta}$ is not integral over~$\BB_H^{\alpha}$.
The question has a positive answer in the following case.

\begin{Prop}\label{integralthm}
If $H= k[G]$ is the Hopf algebra  of a finite group~$G$, then
$S(t_H)_{\Theta}$ is integral over~$\BB_H^{\alpha}$.
\end{Prop}

\pf
It follows from the definition of the comultiplication in~$k[G]$
and the definitions of~$t^{-1}$ and of~$\sigma$ that
\begin{equation}\label{sigma-G}
\sigma(g,h) =
\alpha(g,h)\, \frac{t_g t_h}{t_{gh}} \in \BB_H^{\alpha}
\end{equation}
for all $g, h\in G$.
Since $\alpha(g,h)$ is nonzero for all $g, h\in G$,
the fraction $t_g t_h/t_{gh}$ belongs to~$\BB_H^{\alpha}$.
By an easy induction,
for each grouplike element~$g$ and each integer $k\geq 1$,
we conclude that
${t_g^k}/{t_{g^k}}$ belongs to~$\BB_H^{\alpha}$.
Since by Lagrange's theorem there is $N$ such that $g^N = 1$,
it follows that ${t_g^N}/{t_1}$ belongs to~$\BB_H^{\alpha}$.
Now, $t_1  = \sigma(1,1) \in \BB_H^{\alpha}$ by Lemma~\ref{sigma(1,h)}.
Therefore, ${t_g^N} \in \BB_H^{\alpha}$.
Since by Example~\ref{exempleS(k[G])} of Appendix~\ref{free-comm-coalg},
the elements $t^{\pm 1}_g$ ($g\in G$) generate~$S(t_H)_{\Theta}$,
the conclusion follows.
\epf

\section{The universal twisted algebra}\label{universal twisted algebra}

We have observed in \S~\ref{universal cocycle2} that
the universal cocycle~$\sigma$ takes
values in~$\BB_H^{\alpha}$. 
Therefore we may apply the construction of \S~\ref{Twisted products}
and consider the twisted algebra
\begin{equation}\label{univ-twisted-alg}
\AA_H^{\alpha} = \BB_H^{\alpha} \otimes {}^{\sigma} H \, .
\end{equation}
The product~$\cdot_{\sigma}$ of~$\AA_H^{\alpha}$ is given
for all $b$, $c\in \BB_H^{\alpha}$ and $x$, $y \in H$ by
\begin{equation}\label{twisted-multiplication2}
(b\otimes x) \cdot_{\sigma} (c\otimes y)
= \sum_{(x), (y)}\, bc \, \sigma(x\sw1, y\sw1) \otimes  x\sw2 y\sw2 \, .
\end{equation}

\begin{Lem}\label{sigma-twisted}
The twisted algebra $\AA_H^{\alpha}$ 
is an asso\-ciative unital algebra with unit~$t^{-1}_1 = 1/t_1$.
\end{Lem}

\pf
The associativity follows from the fact that $\sigma$ is a
two-cocycle. The fact that $t^{-1}_1$ is
the unit follows from~\eqref{twisted-multiplication2} and
Lemma~\ref{sigma(1,h)}.
\epf

We call~$\AA_H^{\alpha}$ the \emph{universal twisted algebra}
attached to the two-cocycle~$\alpha$.
Note that the subalgebra of coinvariants of~$\AA_H^{\alpha}$
is equal to~$\BB_H^{\alpha} \otimes k1$; 
this subalgebra is central in~$\AA_H^{\alpha}$.

Consider the restriction $\epsilon : \BB_H^{\alpha} \to k$
of the algebra morphism $\epsilon : S(t_H)_{\Theta} \to k$ of
Lemma~\ref{augmentation} and the maximal ideal 
$$\mm_{\epsilon} = \Ker (\epsilon : \BB_H^{\alpha} \to k) \, .$$
We can recover the original twisted algebra~${}^{\alpha} H$ 
from~$\AA_H^{\alpha}$ as follows.

\begin{Prop}\label{quotient-of-BB}
There is an isomorphism of~$H$-comodule algebras
$$\AA_H^{\alpha}/ \mm_{\epsilon} \, \AA_H^{\alpha} \cong  {}^{\alpha} H \, .$$
\end{Prop}

\pf
We have the obvious isomorphisms of $H$-comodules
$$\AA_H^{\alpha}/ \mm_{\epsilon} \, \AA_H^{\alpha} 
\cong (\BB_H^{\alpha}/ \mm_{\epsilon} ) \otimes H \cong H\, .$$
The fact that the product~$\cdot_{\sigma}$ on~$\AA_H^{\alpha}$
given by~\eqref{twisted-multiplication2} turns into the 
product~$\cdot_{\alpha}$ on~$H$ given by~\eqref{twisted-multiplication}
follows from the equality $\epsilon \circ \sigma = \alpha$ of Lemma~\ref{augmentation}.
\epf

We record the following structure theorem for $\AA_H^{\alpha}$. 
It involves the field of fractions $\KK_H^{\alpha}$ of~$\BB_H^{\alpha}$.

\begin{Thm}\label{AAsemisimple}
Let $k$ be a field of characteristic zero and
$H$ a finite-dimen\-sional Hopf algebra.
If the algebra ${}^{\alpha} H$ is semisimple,
then so is 
$\KK_H^{\alpha}\otimes_{\BB_H^{\alpha}} \AA_H^{\alpha}$.
\end{Thm}

By~\cite[Th.~7.4.2]{M2}, 
if~$H$ is a finite-dimen\-sional semisimple algebra,
then ${}^{\alpha} H$ is semisimple for any two-cocycle~$\alpha$.
Note that there are non-semi\-simple Hopf algebras 
with two-cocycles~$\alpha$ such that ${}^{\alpha} H$ is (semi)\-simple,
e.g., Sweedler's algebra detailed in Section~\ref{Sweedler}.

\pf
Let $A$ be a $B$-algebra, where $B$ is a commutative algebra. 
Suppose that $A$ is free of finite rank
as a $B$-module with basis $\{a_i\}_{i= 1, \ldots , d}$. 
The trace form of~$A$ is the $B$-bilinear form
defined for all $a$, $a' \in A$ by
$$\langle a,a' \rangle_A = \Tr (R_{aa'}) \, ,$$
where $R_{aa'}$ is the right multiplication by~$aa'$ 
and $\Tr : \End_B(A) \to B$ is the trace map. 
It is well known that, if $B$ is a field, then $A$ is semisimple
if and only if the trace form $\langle \, ,\; \rangle_A$ is nondegenerate.
The latter is equivalent to the nonvanishing of the determinant
$\det((\langle a_i,a_j \rangle_A)_{i,j= 1, \ldots , d})$.

Let $f: B \to B'$ be a morphism of commutative algebras.
We denote by the same letter~$f$ the induced algebra morphism $A \to A' = B' \otimes_B A$.
The trace forms of~$A$ and $A'$ are related by the formula
$$\langle f(a), f(a') \rangle_{A'} = f(\langle a,a' \rangle_A)$$
for all $a$, $a' \in A$.

We apply this to the case where $B = \BB_H^{\alpha}$ and $A = \AA_H^{\alpha}$.
If $\{h_i\}_{i= 1, \ldots , d}$ is a basis of~$H$, then it is also a basis of~$\AA_H^{\alpha}$
considered as a $\BB_H^{\alpha}$-module.
Set 
$$D = \det((\langle h_i,h_j \rangle_{\AA_H^{\alpha}})_{i,j= 1, \ldots , d}) \in \BB_H^{\alpha}\, .$$
By the remarks above, the trace form of $\KK_H^{\alpha}\otimes_{\BB_H^{\alpha}} \AA_H^{\alpha}$
is the image of the trace form of~$\AA_H^{\alpha}$ in~$\KK_H^{\alpha}$.
Therefore, $\KK_H^{\alpha}\otimes_{\BB_H^{\alpha}} \AA_H^{\alpha}$ is
semisimple if~$D \neq 0$ in~$\KK_H^{\alpha}$, hence in~$\BB_H^{\alpha}$. 
In order to prove the nonvanishing of~$D$,
it is enough to map it to a nonzero element.
Consider the augmentation map $\epsilon : \BB_H^{\alpha} \to k$ of 
\S~\ref{universal cocycle2}.
By Proposition~\ref{quotient-of-BB}, 
$$k \otimes_{\BB_H^{\alpha}} \AA_H^{\alpha} 
= \AA_H^{\alpha}/ \mm_{\epsilon} \,  \AA_H^{\alpha} \cong {}^{\alpha} H \, .$$
The trace form of~${}^{\alpha} H$ is the image under~$\epsilon$ 
of the trace form of~$\AA_H^{\alpha}$.
Now,
$$\epsilon(D) = \det((\langle h_i,h_j \rangle_{{}^{\alpha} H})_{i,j= 1, \ldots , d}) \in k $$
is nonzero since ${}^{\alpha} H$ is semisimple. Therefore,~$D \neq 0$.
\epf

Below we shall show under the same hypotheses as in Theorem~\ref{AAsemisimple}
that $\KK_H^{\alpha}\otimes_{\BB_H^{\alpha}} \AA_H^{\alpha}$
is a simple algebra if ${}^{\alpha} H$ is simple.
Let us first establish the following result.

\begin{Lem}\label{centerAA}
If the center of the algebra ${}^{\alpha} H$ is trivial, that is, 
consists of the scalar multiples of the unit, 
then the center of~$\AA_H^{\alpha}$ is equal to~$\BB_H^{\alpha}\otimes 1$.
\end{Lem}

\pf
Since the universal cocycle~$\sigma$ takes values in~$S(t_H)_{\Theta}$ 
(see \S~\ref{universal cocycle2}),
we may apply the construction of \S~\ref{Twisted products}
and consider the twisted algebra
\begin{equation}\label{univ-twisted-alg2}
S(t_H)_{\Theta} \otimes {}^{\sigma} H \, ,
\end{equation}
which is clearly isomorphic to the algebra
$S(t_H)_{\Theta} \otimes_{\BB_H^{\alpha}} \AA_H^{\alpha}$.
It follows from Lemma~\ref{sigma-twisted} that $S(t_H)_{\Theta} \otimes {}^{\sigma} H$ 
is an asso\-ciative algebra with unit~$t^{-1}_1 = 1/t_1$.

Equation~\eqref{sigma-def} makes sense in~$S(t_H)_{\Theta}$.
It can be interpreted as saying that the two-cocycles~$\sigma$ 
and~$\alpha$ are cohomologous in~$S(t_H)_{\Theta}$.
By~\cite{BCM} and~\cite{Doi} (see also~\cite[Sect.~7]{M2}), the linear map 
$$\varphi : S(t_H)_{\Theta} \otimes {}^{\sigma} H 
\to S(t_H)_{\Theta} \otimes {}^{\alpha} H$$
given for all $b\in S(t_H)_{\Theta}$ and $h\in H$ by
\begin{equation}\label{formula-iso}
\varphi(b\otimes h) = \sum_{(h)}\, b\, t_{h_{(1)}} \otimes h_{(2)}
\end{equation}
is a $S(t_H)_{\Theta}$-linear isomorphism of $H$-comodule algebras.
Note that the inverse map 
$$\varphi^{-1} : S(t_H)_{\Theta} \otimes {}^{\alpha} H \to S(t_H)_{\Theta} \otimes {}^{\sigma} H$$
is given for all $b\in S(t_H)_{\Theta}$ and $h\in H$ by
\begin{equation}\label{formula-inv-iso}
\varphi^{-1}(b\otimes h) = \sum_{(h)}\, b\, t^{-1}_{h_{(1)}} \otimes h_{(2)} \, .
\end{equation}
The center of~$S(t_H)_{\Theta} \otimes {}^{\sigma} H$ is isomorphic
to the center of~$S(t_H)_{\Theta} \otimes {}^{\alpha} H$ under~$\varphi$. 
Now it is easy to check that,
under the hypothesis of the lemma, the center of~$S(t_H)_{\Theta} \otimes {}^{\alpha} H$
is the subalgebra~$S(t_H)_{\Theta} \otimes 1$.
The isomorphism~$\varphi$ being $S(t_H)_{\Theta}$-linear and 
sending the unit~$t^{-1}_1 \otimes  1 \in S(t_H)_{\Theta} \otimes {}^{\sigma} H$ to
the unit~$1 \otimes 1\in S(t_H)_{\Theta} \otimes {}^{\alpha} H$,
the center of~$S(t_H)_{\Theta} \otimes {}^{\sigma} H$ is equal to
$$t^{-1}_1S(t_H)_{\Theta} \otimes 1 = S(t_H)_{\Theta} \otimes 1\, .$$

Now, let $\omega \in \AA_H^{\alpha} = \BB_H^{\alpha} \otimes {}^{\sigma} H$ be a central element. 
It remains central when considered in~$S(t_H)_{\Theta} \otimes {}^{\sigma} H$.
By the above argument, we have $\omega = b\otimes 1$ for some $b\in S(t_H)_{\Theta}$.
Since $\BB_H^{\alpha} \subset S(t_H)_{\Theta}$, 
it follows that~$b\in\BB_H^{\alpha}$.
\epf

\begin{Cor}\label{AAsimple}
Let $k$ be a field of characteristic zero and $H$ a finite-dimen\-sional Hopf algebra.
If the algebra ${}^{\alpha} H$ is simple,
then so is the algebra $\KK_H^{\alpha}\otimes_{\BB_H^{\alpha}} \AA_H^{\alpha}$.
\end{Cor}

\pf
By Theorem~\ref{AAsemisimple}, 
$\KK_H^{\alpha}\otimes_{\BB_H^{\alpha}} \AA_H^{\alpha}$ is semisimple.
By Lemma~\ref{centerAA}, the center of~$\AA_H^{\alpha}$
is~$\BB_H^{\alpha} \otimes 1$. Therefore,
the center of~$K_H \otimes_{\BB_H^{\alpha}} \AA_H^{\alpha}$
is~$K_H \otimes 1$. 
The conclusion follows immediately.
\epf

Under the hypotheses of Corollary~\ref{AAsimple}, we may wonder
what the \emph{index} of the central simple algebra 
$\KK_H^{\alpha}\otimes_{\BB_H^{\alpha}} \AA_H^{\alpha}$
is and how it depends on the two-cocycle~$\alpha$.
Such questions are open even in the case of group algebras; 
for partial answers, see~\cite[\S~3]{AHN1}.

\section{Forms of~${}^{\alpha} H$}\label{forms}

As in Section~\ref{universal cocycle}, we fix a Hopf algebra~$H$
and a normalized convolution-invertible two-cocycle $\alpha : H \times H \to k$.
We consider the corresponding twisted algebra~${}^{\alpha} H$.

Given a normalized convolution-invertible two-cocycle $\beta : H \times H \to K$
with values in a (field) extension $K$ of~$k$,
we say that the twisted $H$-comodule algebra $K \otimes {}^{\beta} H$
is a \emph{$K$-form} of ${}^{\alpha} H$ if
there is an extension~$L$ of~$K$ and 
an $L$-linear isomorphism of~$H$-comodule algebras
$$L \otimes_K (K \otimes {}^{\beta} H) \cong L\otimes_k {}^{\alpha} H \, .$$

We now state two theorems relating forms of~${}^{\alpha} H$
to central specializations of
the universal twisted algebra~$ \AA_H^{\alpha}$ introduced in the previous section.

\begin{Thm}\label{formthm2}
For any $K$-form $K \otimes {}^{\beta} H$ of ${}^{\alpha} H$,
where $\beta : H \times H \to K$ is a normalized convolution-invertible two-cocycle
with values in an extension~$K$ of~$k$,
there are an algebra morphism $\lambda : \BB_H^{\alpha} \to K$ 
and a $K$-linear isomorphism of~$H$-comodule algebras
$$K \otimes_{\BB_H^{\alpha}} \AA_H^{\alpha} \cong K \otimes {}^{\beta} H \, .$$
\end{Thm}

The algebra morphism $\lambda : \BB_H^{\alpha} \to K$ equips $K$ with
a natural $\BB_H^{\alpha}$-module structure so that it makes sense
to consider the algebra $K \otimes_{\BB_H^{\alpha}} \AA_H^{\alpha}$
as in the previous statement.

\pf
Since there is an $L$-linear isomorphism of~$H$-comodule algebras
$$L \otimes_K (K \otimes {}^{\beta} H) \cong L\otimes_k {}^{\alpha} H $$
for some extension~$L$ of~$K$,
by~\cite{BCM},~\cite{DT2} 
there is a convolution invertible linear map
$\lambda_0 : H \to L$ such that for all $g$, $h\in H$,
\begin{equation}\label{beta-alpha}
\beta(x,y) = \sum_{(x), (y)}\, \lambda_0(x_{(1)}) \,
\lambda_0(y_{(1)}) \, \alpha(x_{(2)},y_{(2)}) \,
\lambda_0^{-1}(x_{(3)} y_{(3)}) \, ,
\end{equation}
where $\lambda_0^{-1}$ is the convolution inverse of~$\lambda_0$. 
Define an algebra morphism $\lambda : S(t_H) \to L$ 
by $\lambda(t_h) = \lambda_0(h)$ for all $h\in H$. 
This morphism extends to an algebra morphism, which we still denote
by~$\lambda$, from $S(t_H)_{\Theta}$ to~$L$ such that
$\lambda(t^{-1}_h) = \lambda_0^{-1}(h)$ for all $h\in H$.
We obtain an algebra morphism $\lambda : \BB_H^{\alpha} \to L$ by
restricting~$\lambda$ to the subalgebra~$\BB_H^{\alpha}$.

It follows from \eqref{sigma-def}, \eqref{sigma^{-1}-def}, and \eqref{beta-alpha} that
$\lambda$ sends the generators $\sigma(x,y)$ and $\sigma^{-1}(x,y)$
of~$\BB_H^{\alpha}$ respectively to~$\beta(x,y)$ and to~$\beta^{-1}(x,y)$,
where $\beta^{-1}$ is the convolution inverse of~$\beta$.
Since $\beta(x,y)$ and $\beta^{-1}(x,y)$ belong to~$K$,
we conclude that
$\lambda(\BB_H^{\alpha}) \subset K$.
Thus, $\lambda\otimes \id_H$ defines a map
$$\AA_H^{\alpha} = \BB_H^{\alpha} \otimes {}^{\sigma} H \to K \otimes {}^{\beta} H \, .$$
Since $\lambda \circ \sigma = \beta$, this map is an algebra morphism;
it is also a comodule map. 
Therefore, the algebra morphism $\lambda: \BB_H^{\alpha} \to K$
induces an $H$-comodule algebra isomorphism
$K \otimes_{\BB_H^{\alpha}} \AA_H^{\alpha} \cong K \otimes {}^{\beta} H \, .$
\epf

Under the integrality condition of~\S~\ref{integrality},
Theorem~\ref{formthm2} has the following converse.

\begin{Thm}\label{formthm}
If the algebra $S(t_H)_{\Theta}$ is integral over the subalgebra~$\BB_H^{\alpha}$, 
then for any extension $K$ of~$k$ and any algebra morphism
$\lambda : \BB_H^{\alpha} \to K$, the $H$-comodule $K$-algebra
$K \otimes_{\BB_H^{\alpha}} \AA_H^{\alpha}$
is a $K$-form of~${}^{\alpha} H$.
\end{Thm}

\pf
Let $\mm$ be the kernel of an algebra morphism $\lambda : \BB_H^{\alpha} \to K$;
it is a maximal ideal of~$\BB_H^{\alpha}$,
and the image~$k' = \BB_H^{\alpha}/\mm$ of~$\lambda$ is an extension of~$k$ 
contained in~$K$.
It follows from the integrality assumption
that there is a maximal ideal~$\MM$ of~$S(t_H)_{\Theta}$
such that $\mm \, S(t_H)_{\Theta} \subset \MM$ 
(see e.g., \cite[Chap.~IX, Proposition~9]{La}).
Let $\widetilde{\lambda}$ be the projection $S(t_H)_{\Theta} \to L = S(t_H)_{\Theta}/\MM$.
We have the following commutative square in which the top horizontal map is the 
inclusion $\BB_H^{\alpha} \hookrightarrow S(t_H)_{\Theta}$:
$$\begin{CD}
\BB_H^{\alpha} @>{}>> & S(t_H)_{\Theta}  \\
@V{\lambda}VV & @VV{\widetilde{\lambda}}V  \\
k' @>>> & L  \\
\end{CD}$$
The bottom horizontal map $k' \to L$ in the previous square is injective
since $k'$ and $L$ are fields.

Applying $\widetilde{\lambda}$ to \eqref{sigma-def} and setting
$\beta = \widetilde{\lambda} \circ \sigma : H\times H \to L$, we obtain
\begin{eqnarray*}
\beta(g,h) & = &
\widetilde{\lambda}\bigl(\sigma(g,h)\bigr) \\
& = &\sum_{(g), (h)}\, \widetilde{\lambda}\bigl(t_{g_{(1)}} \bigr)
\widetilde{\lambda} \bigl( t_{h_{(1)}} \bigr) \, \widetilde{\lambda}
\bigl(\alpha(g_{(2)},h_{(2)}) \bigr) \, \widetilde{\lambda} \bigl(
t^{-1}_{g_{(3)} h_{(3)}} \bigr)
\end{eqnarray*}
for all $g$, $h\in H$. Since $\widetilde{\lambda}$ is an algebra
morphism, $\widetilde{\lambda} \circ t^{-1}$ is the convolution inverse
of $\widetilde{\lambda} \circ t$. 
It follows that $\beta$ is an invertible two-cocycle, which is 
cohomologous to~ $\widetilde{\lambda} \circ \alpha = \alpha$. 
Therefore,
$$L \otimes {}^{\alpha} H \cong L \otimes {}^{\beta} H 
= L\otimes_K \bigl(K \otimes_{\BB_H^{\alpha}} (\BB_H^{\alpha} \otimes {}^{\sigma} H) \bigr)
= L \otimes_{\BB_H^{\alpha}} \AA_H^{\alpha}\, , $$
which shows that $K \otimes_{\BB_H^{\alpha}} \AA_H^{\alpha}$
is a $K$-form of~$ {}^{\alpha} H$. 
\epf

\begin{Cor}\label{formcoro}
If $S(t_H)_{\Theta}$ is integral over~$\BB_H^{\alpha}$,
then for any extension $K$ of~$k$ there is a surjective map
$$\Alg(\BB_H^{\alpha},K) \to \Forms_{K}({}^{\alpha} H)$$
from the set of algebra morphisms $\BB_H^{\alpha} \to K$ to the set
of isomorphism classes of $K$-forms of~${}^{\alpha} H$.
\end{Cor}

\pf 
For any algebra morphism $\lambda : \BB_H^{\alpha} \to K$, the
$H$-comodule $K$-algebra $\lambda_*(\BB_H^{\alpha}) =
K \otimes_{\BB_H^{\alpha}} \AA_H^{\alpha}$ is a $K$-form
of~${}^{\alpha} H$ by Theorem~\ref{formthm}. 
By Theorem~\ref{formthm2},
the map $\lambda \mapsto \lambda_*(\BB_H^{\alpha})$ maps $\Alg(\BB_H^{\alpha},K)$ 
surjectively onto the set~$\Forms_{K}({}^{\alpha} H)$. 
\epf

This corollary allows us to consider the cleft $H$-Galois extension
$\BB_H^{\alpha} \subset \AA_H^{\alpha}$ as a
\emph{versal deformation space} for the forms of~${}^{\alpha} H$
(this space is not universal since the algebra morphism
$\lambda : \BB_H^{\alpha} \to K$ in Theorem~\ref{formthm2}
is not necessarily unique).
As a further consequence, if
$\BB_H^{\alpha} \cong k[u_1, \ldots, u_r]/(P_1, \ldots, P_s)$,
then the polynomials $P_1, \ldots, P_s \in k[u_1, \ldots, u_r]$ 
are invariants of forms of~${}^{\alpha} H$.

In case the twisted algebra ${}^{\alpha} H$ is simple,
we have another important consequence of Theorem~\ref{formthm}.

\begin{Thm}\label{Azumaya}
If $S(t_H)_{\Theta}$ is integral over~$\BB_H^{\alpha}$
and if the algebra ${}^{\alpha} H$ is simple,
then $\AA_H^{\alpha}$ is an Azumaya algebra.
\end{Thm}

\pf
By Lemma~\ref{centerAA},
$\BB_H^{\alpha}$ is the center of~$\AA_H^{\alpha}$.
By~\cite[Th.~15]{Az}, in order to prove the theorem, it suffices to
check that $\AA_H^{\alpha}/\mm \, \AA_H^{\alpha}$ is simple
for any maximal ideal~$\mm$ of~$\BB_H^{\alpha}$.
Now, by Theorem~\ref{formthm},
$$\AA_H^{\alpha}/\mm \, \AA_H^{\alpha} 
= (\BB_H^{\alpha}/\mm) \otimes_{\BB_H^{\alpha}} \AA_H^{\alpha}$$
is a $\BB_H^{\alpha}/\mm$-form of the simple algebra~${}^{\alpha} H$. Therefore, it is simple.
\epf

\section{Relating the universal comodule algebra 
to the universal twisted algebra}\label{forms-identities}

From now on we assume that the ground field is infinite.
As in Section~\ref{universal cocycle}, we fix a Hopf algebra~$H$
and a normalized convolution-invertible two-cocycle $\alpha : H \times H \to k$,
and we consider the corresponding twisted algebra~${}^{\alpha} H$.
Since~${}^{\alpha} H$ satisfies the assumptions of \S~\ref{assumptions}
with $B=k$ and~$u = \id_{H}$,
we may apply the constructions and 
the results of Section~\ref{universal evaluation} to~${}^{\alpha} H$.

Let  
$\mu_{\alpha}: T(X_H) \to S(t_H) \otimes {}^{\alpha} H$
be the universal evaluation morphism defined by~\eqref{universal map}.
In the present situation,
\begin{equation}\label{universal map2}
\mu_{\alpha}(X_h) = \sum_{(h)}\, t_{h\sw1} \otimes h\sw2
\end{equation}
for all $h \in H$. 
For short, we denote~$ \UU_H({}^{\alpha} H)$ by~$\UU_H^{\alpha}$.
By Corollary~\ref{inject},
$\mu_{\alpha}$~induces an injective $H$-comodule algebra morphism
$$\UU_H^{\alpha} \hookrightarrow S(t_H) \otimes {}^{\alpha} H \, , $$
which we still denote~$\mu_{\alpha}$.
It follows from Proposition~\ref{coinv-UU}
that $P \in \UU_H^{\alpha}$ is coinvariant if and only if 
$\mu_{\alpha}(P) \in S(t_H) \otimes 1$.

Recall the $H$-comodule algebra isomorphism
$$\varphi^{-1} : S(t_H)_{\Theta} \otimes {}^{\alpha} H \to S(t_H)_{\Theta} \otimes {}^{\sigma} H$$
given by~\eqref{formula-inv-iso}.
Define $\mu_{\sigma} : \UU_H^{\alpha} \to S(t_H)_{\Theta} \otimes {}^{\sigma} H$ by
\begin{equation}\label{def-mu-sigma}
\mu_{\sigma} = \varphi^{-1} \circ \mu_{\alpha}\, ,
\end{equation}
where we now view $\mu_{\alpha}$ as a morphism
with values in~$S(t_H)_{\Theta} \otimes  {}^{\alpha} H$.
It follows from the definition and from Corollary~\ref{inject}
that $\mu_{\sigma}$
is an injective $H$-comodule algebra morphism.

\begin{Lem}\label{lem-mu-sigma}
(a) We have $\mu_{\sigma}(X_h) = 1 \otimes h$ for all $h\in H$.

(b) The morphism $\mu_{\sigma}$
maps $\UU_H^{\alpha}$ into~$\AA_H^{\alpha}$.
\end{Lem}

\pf (a) It suffices to check that
$\varphi(1\otimes h)  = \mu_{\alpha}(X_h)$
for all $h\in H$. 
These equalities follow from~\eqref{formula-iso}
and~\eqref{universal map2}. 

(b) By Part\,(a) and by \eqref{twisted-multiplication2},
$$\mu_{\sigma}(X_g X_h)
=  \sum_{(g), (h)}\, \sigma(g\sw1, h\sw1) \otimes g\sw2 h\sw2 $$
for all $g$, $h\in H$.
Using~\eqref{twisted-multiplication2} repeatedly, we see that
$\mu_{\sigma}$ sends all monomials in the variables~$X_h$,
hence all elements of~$\UU_H^{\alpha}$,
to~$\BB_H^{\alpha} \otimes {}^{\sigma} H = \AA_H^{\alpha}$.
\epf

Let $Z({}^{\alpha} H)$ be the center of~${}^{\alpha} H$.
We characterize the center~$\ZZ_H^{\alpha}$ of $\UU_H^{\alpha}$
as follows.

\begin{Prop}\label{mu-alpha-center}
An element $w\in\UU_H^{\alpha}$ belongs to the center~$\ZZ_H^{\alpha}$
if and only if~$\mu_{\alpha}(w)$ belongs to~$S(t_H)\otimes Z({}^{\alpha} H)$.
\end{Prop}

\pf
Suppose first that $\mu_{\alpha}(w) \in S(t_H)\otimes Z({}^{\alpha} H)$.
The elements of the latter are clearly central in~$S(t_H) \otimes {}^{\alpha} H$.
Therefore, $\mu_{\alpha}(w)$ commutes with all elements~$\mu_{\alpha}(X_h)$, where
$h$~runs over~$H$. 
Since $\mu_{\alpha}$ is an injective algebra morphism, $w$ commutes with all generators~$X_h$ 
of~$\UU_H^{\alpha}$. Hence, $\mu_{\alpha}(w) \in \ZZ_H^{\alpha}$.

Conversely, any element $w \in \ZZ_H^{\alpha}$
commutes with all generators~$X_h$ of~$\UU_H^{\alpha}$.
Therefore, by Lemma~\ref{lem-mu-sigma}, 
its image $\mu_{\sigma}(w)$ commutes with $\mu_{\sigma}(X_h) = 1 \otimes h$ for all $h\in H$.
Since $S(t_H)_{\Theta} \otimes {}^{\sigma}H$ is generated as an algebra by~$1 \otimes H$
and the central subalgebra~$S(t_H)_{\Theta} \otimes 1$, it follows immediately
that $\mu_{\sigma}(w)$ commutes with all elements of~$S(t_H)_{\Theta} \otimes {}^{\sigma}H$.
Consequently, $\mu_{\alpha}(w) = \varphi(\mu_{\sigma}(w))$ commutes
with all elements of~$S(t_H)_{\Theta} \otimes {}^{\alpha}H$,
i.e., belongs to the center of the latter.
Write 
$$\mu_{\alpha}(w) = \sum_{j\in J}\, s_j \otimes h_j\, ,$$
where $\{ s_j \}_{j\in J}$ is a basis of~$S(t_H)_{\Theta}$ over the ground field
and $h_j \in H$ for all $j\in J$.
Using brackets to denote commutators, we obtain
$$\sum_{j\in J}\, s_j \otimes [h_j,h] = \Bigl[\sum_{j\in J}\, s_j \otimes h_j,h \Bigr] = 
[\mu_{\alpha}(w), 1 \otimes h]$$
for all $h\in H$.
This implies that $[h_j,h] = 0$ for all $h\in H$ and~$j\in J$.
Therefore, $h_j \in Z({}^{\alpha} H)$ for all~$j\in J$,
from which it follows that 
$\mu_{\alpha}(w)$ lies in~$S(t_H)_{\Theta} \otimes Z({}^{\alpha} H)$.
Since $\mu_{\alpha}(w) \in S(t_H) \otimes {}^{\alpha} H$,
it lies in~$S(t_H) \otimes Z({}^{\alpha} H)$.
\epf

\begin{Cor}\label{Znozerodivisors}
The center $\ZZ_H^{\alpha}$ of~$\UU_H^{\alpha}$
has no zero-divisors if~$Z({}^{\alpha} H)$ has none.
\end{Cor}

\pf
By Corollary~\ref{inject} and Proposition~\ref{mu-alpha-center}, 
$\mu_{\alpha}$ embeds~$\ZZ_H^{\alpha}$ into
the algebra~$S(t_H)_{\Theta} \otimes Z({}^{\alpha} H)$, 
which has no zero-divisors if~$Z({}^{\alpha} H)$ has none.
\epf

\section{The case of nondegenerate cocycles}\label{nondegenerate}

Let $H$ be a Hopf algebra and $\alpha : H \times H \to k$
a normalized convolution-invertible two-cocycle.
Generalizing terminology used for group cocycles, 
we say that $\alpha$ is \emph{nondegenerate} 
if the center of the twisted algebra~${}^{\alpha} H$ is trivial,
i.e., if $Z({}^{\alpha} H) = k1$.

Restricting to nondegenerate two-cocycles allows us 
to relate~$\ZZ_H^{\alpha}$ to $\BB_H^{\alpha}$ as follows.

\begin{Prop}\label{mu-sigma-center}
If $\alpha$ is a nondegenerate two-cocycle,
then $\mu_{\alpha}$ and $\mu_{\sigma}$
map $\ZZ_H^{\alpha}$ into~$\BB_H^{\alpha} \otimes 1$.
\end{Prop}

\pf
Since $\alpha$ is nondegenerate, it follows from Proposition~\ref{mu-alpha-center}
that $\mu_{\alpha}$ maps $\ZZ_H^{\alpha}$ into
$S(t_H)\otimes Z({}^{\alpha} H) = S(t_H)\otimes 1$.
The isomorphism 
$$\varphi^{-1} : S(t_H)_{\Theta} \otimes {}^{\alpha} H \to S(t_H)_{\Theta} \otimes {}^{\sigma} H$$
being $S(t_H)_{\Theta}$-linear,  
$\mu_{\sigma} = \varphi^{-1} \circ \mu_{\alpha}$ maps $\ZZ_H^{\alpha}$ into 
$$\varphi^{-1}(S(t_H)\otimes 1) = S(t_H)\, \varphi^{-1}(1\otimes 1) = t^{-1}_1 S(t_H)\otimes 1 \, ,$$
which is contained in~$S(t_H)_{\Theta} \otimes 1$.
Now, by Lemma~\ref{lem-mu-sigma}\,(b), $\mu_{\sigma}$
maps $\UU_H^{\alpha}$ into~$\AA_H^{\alpha} = \BB_H^{\alpha} \otimes {}^{\sigma} H$.
Since $\BB_H^{\alpha} \subset S(t_H)_{\Theta}$, we conclude that
$$\mu_{\sigma}(\ZZ_H^{\alpha}) \subset \BB_H^{\alpha} \otimes 1\, .$$
To obtain the same inclusion with $\mu_{\alpha}$ in lieu of~$\mu_{\sigma}$,
we apply the $S(t_H)_{\Theta}$-linear isomorphism~$\varphi$ to the previous
inclusion. We obtain
$$\mu_{\alpha}(\ZZ_H^{\alpha}) \subset \BB_H^{\alpha} \, \varphi(1\otimes 1)
= t_1 \BB_H^{\alpha} \otimes 1 \, .$$
We conclude the proof by observing that $t_1$ is an invertible element of~$\BB_H^{\alpha}$
as a consequence of Lemma~\ref{sigma(1,h)} applied to~$x=1$.
\epf

If $\alpha$ is a nondegenerate two-cocycle,
as a consequence of Proposition~\ref{mu-sigma-center}, we can view the center
$\ZZ_H^{\alpha}$ of~$\UU_H^{\alpha}$ as a subalgebra of~$\BB_H^{\alpha}$,
and consider the algebra
\begin{equation}
\BB_H^{\alpha} \otimes_{\ZZ_H^{\alpha}} \UU_H^{\alpha} \, .
\end{equation}

\begin{Rem}
In~\cite{AEGN} nondegenerate cocycles are defined in a more restrictive way than above, 
namely as those cocycles~$\alpha$ such that the algebra~$^{\alpha}H$ is $k$-central simple. 
Both definitions coincide when $H$ is semisimple. 
Furthermore, one can show that, when $H=k[G]$ is the Hopf algebra of a group, 
there exist nondegenerate two-cocycles only if $H$ is semisimple.
\end{Rem}

Our next objective is to relate 
$\BB_H^{\alpha} \otimes_{\ZZ_H^{\alpha}} \UU_H^{\alpha}$
to the universal twisted algebra~$\AA_H^{\alpha}$.
This will be possible under the following additional condition on the two-cocycle:
a nondegenerate two-cocycle~$\alpha :  H \times H \to k$
is said to be \emph{nice} 
if $\BB_H^{\alpha}$ is a localization of~$\mu_{\alpha} (\ZZ_H^{\alpha})$
or, equivalently, if their fields of fractions coincide.
If $\alpha$ is nice, then $\BB_H^{\alpha}$ is flat as a $\ZZ_H^{\alpha}$-module
and 
\begin{equation}\label{BBtBB}
\BB_H^{\alpha} \otimes_{\ZZ_H^{\alpha}} \BB_H^{\alpha} \cong \BB_H^{\alpha} \, .
\end{equation}

\begin{Thm}\label{UGalois}
If $\alpha : H\times H \to k$ is a nice nondegenerate two-cocycle,
then $\mu_{\sigma}$ induces an isomorphism of $H$-comodule algebras
$$\BB_H^{\alpha} \otimes_{\ZZ_H^{\alpha}} \UU_H^{\alpha} \cong \AA_H^{\alpha} \, .$$
\end{Thm}

\pf
The $H$-comodule algebra morphism
$\mu_{\sigma}: \UU_H^{\alpha} \to \BB_H^{\alpha} \otimes {}^{\sigma} H$
induces an $H$-comodule algebra morphism
$$ \id\otimes \mu_{\sigma} : \BB_H^{\alpha} \otimes_{\ZZ_H^{\alpha}}\UU_H^{\alpha} 
\to \BB_H^{\alpha} \otimes_{\ZZ_H^{\alpha}}\BB_H^{\alpha} \otimes {}^{\sigma} H \, .$$
Since $\mu_{\sigma}$ is injective and $\BB_H^{\alpha}$ is flat as a $\ZZ_H^{\alpha}$-module,
$\id\otimes \mu_{\sigma}$ is injective. 
It then follows from~\eqref{BBtBB} that $\id\otimes \mu_{\sigma}$ maps 
$\BB_H^{\alpha} \otimes_{\ZZ_H^{\alpha}}\UU_H^{\alpha} $ injectively into
$\BB_H^{\alpha} \otimes {}^{\sigma} H = \AA_H^{\alpha}$.

As an $\BB_H^{\alpha}$-algebra, $\AA_H^{\alpha}$ is generated by
the elements~$1\otimes h$, where $h$ runs over~$H$.
By Lemma~\ref{lem-mu-sigma}\,(a), these elements are in the image of~$\mu_{\sigma}$.
This proves the surjectivity of~$\id\otimes \mu_{\sigma}$.
\epf

The following is a consequence of Lemma~\ref{centerAA} and Theorem~\ref{UGalois}.

\begin{Cor}
If $\alpha$ is a nice nondegenerate two-cocycle,
then the center 
of~$\BB_H^{\alpha} \otimes_{\ZZ_H^{\alpha}} \UU_H^{\alpha}$ is~$\BB_H^{\alpha}$.
\end{Cor}

We do not know if all nondegenerate two-cocycles are nice.
In the group case, we have the following result.

\begin{Prop}\label{G-nice}
Any two-cocycle on the Hopf algebra $k[G]$ of a group is nice.
\end{Prop}

\pf
For $g$, $h\in G$, set $Z_g = X_g X_{g^{-1}}$
and $Z_{g,h} = X_g X_h X_{h^{-1}g^{-1}} \in T(X_H)$, where~$H = k[G]$.
Let $\alpha : H\times H \to k$ be a normalized convolution invertible two-cocycle.
Its restriction to~$G \times G$ takes nonzero values.
An easy computation shows that
\begin{equation}\label{Zg}
\mu_{\alpha}(Z_g)  = \alpha(g,g^{-1})\, t_g \, t_{g^{-1}} \otimes 1
\end{equation}
and
\begin{equation}\label{Zgh}
\mu_{\alpha}(Z_{g,h})  = \alpha(g,g^{-1})\, \alpha(h,(gh)^{-1})\, t_g \, t_h \, t_{(gh)^{-1}} \otimes 1\, .
\end{equation}
It follows from these identities and from Proposition~\ref{mu-alpha-center} that
$Z_g$ and $Z_{g,h}$ represent elements of the center~$\ZZ_H^{\alpha}$.
As an easy consequence of \eqref{sigma-G},~\eqref{Zg}, and~\eqref{Zgh}, 
we obtain the following identities for all~$g$, $h\in G$:
$$\sigma(g,h) = 
\frac{\alpha(g,h)\, \alpha(gh,(gh)^{-1})}{\alpha(g,g^{-1})\, \alpha(h,(gh)^{-1})}
\cdot \frac{\mu_{\alpha}(Z_{g,h})}{\mu_{\alpha}(Z_{gh})}\, .$$
This shows that the generators $\sigma(g,h)$ and $\sigma^{-1}(g,h) = 1/\sigma(g,h)$
of~$\BB_H^{\alpha}$ can be expressed as fractions of elements of~$\mu_{\alpha} (\ZZ_H^{\alpha})$.
\epf

\section{The Sweedler algebra}\label{Sweedler}

We assume in this section that the characteristic of the ground
field~$k$ is different from~2.

\subsection{Definition}

The \emph{Sweedler algebra} $H_4$ is the
algebra generated by two elements $x$, $y$ subject to the relations
\begin{equation}\label{H4def}
x^2 = 1 \, , \quad y^2 = 0 \, , \quad xy + yx = 0 \, .
\end{equation}
It is four-dimensional with basis $\{1, x, y, z\}$, where $z = xy$.

It carries a structure of Hopf algebra with
comultiplication~$\Delta$, counit~$\eps$, and antipode~$S$ given by
\begin{equation}\label{H4comult}
\begin{matrix}
\Delta(1) &=& 1 \otimes 1 \, , & \Delta(x) &=& x \otimes x \, , \\
\Delta(y) &=& 1 \otimes y + y \otimes x \, , & \Delta(z) &=& x \otimes z + z \otimes 1 \, ,\\
\eps(1) &=& \eps(x) = 1 \, , & \eps(y) &=& \eps(z) = 0 \, ,\\
S(1) = 1 \, , && S(x) = x \, , & S(y) = z \, , && S(z) = -y \, .
\end{matrix}
\end{equation}

The tensor algebra $T(X_{H_4})$ is freely generated by the four elements
$$X_1 = E \, , \;\; X_x = X \, , \;\; X_y = Y \, , \;\; X_z  = Z\, ,$$ 
whereas $S(t_{H_4})$ is the
polynomial algebra over the elements $t_1$, $t_x$, $t_y$, $t_z$.

By \eqref{grouplike-tbar} and \eqref{primitive-tbar},
the elements $t^{-1}_1$, $t^{-1}_x$, $t^{-1}_y$, $t^{-1}_z$
of the field of fractions~$K_{H_4}$ of~$S(t_{H_4})$ are given by
\begin{equation*}
t^{-1}_1 = \frac{1}{t_1} \, ,\quad  
t^{-1}_x = \frac{1}{t_x} \, ,\quad  
t^{-1}_y =  - \frac{t_y}{t_1t_x}\ , \quad  
t^{-1}_z = - \frac{t_z} {t_1 t_x}  \, .
\end{equation*}
It follows that the algebra $S(t_{H_4})_{\Theta}$ is isomorphic to
the algebra of Laurent polynomials
$$S(t_{H_4})_{\Theta} \cong k[t_1, t_1^{-1}, t_x, t_x^{-1}, t_y, t_z] \, .$$
(This also follows from Proposition~\ref{Takeuchi}.)

\subsection{Twisted algebras}
Masuoka~\cite{Ma1} classified all Hopf Galois extensions for~$H_4$
(we also follow~\cite{DT2}). In particular, he showed that any
$H_4$-Galois extension of~$k$ is, up to isomorphism,
of the form ${}^{\alpha} H_4$, 
where $\alpha$ is the normalized convolution
invertible two-cocycle given by
\begin{equation}\label{H4alpha}
\begin{matrix}
\alpha(x,x) &=& a \, , \quad \alpha(x,y) &=& 0 \, ,  \quad \alpha(x,z) &=& 0 \, , \\
\alpha(y,x) &=& b \, , \quad \alpha(y,y) &=& c \, , \quad  \alpha(y,z) &=& -c \, , \\
\alpha(z,x) &=& b \, , \quad \alpha(z,y) &=& c \, , \quad  \alpha(z,z)
&=&-ac
\end{matrix}
\end{equation}
for some $a$, $b$, $c \in k$ with $a\neq 0$.
The four-dimensional algebra ${}^{\alpha} H_4$ has the same basis
as~$H_4$ and $1$ is its unit. In this basis the multiplication~$\cdot_{\alpha}$
of~${}^{\alpha} H_4$ is given by
\begin{equation}\label{tablemultiplic}
\begin{matrix}
x \cdot_{\alpha} x &=& a \, , \quad &x \cdot_{\alpha} y &=& z \, , 
\quad &x \cdot_{\alpha} z &=& ay \, , \\
y \cdot_{\alpha} x &=& b - z \, , \quad &y \cdot_{\alpha} y &=& c \, , 
\quad &y \cdot_{\alpha} z &=& -cx + by \, , \\
z \cdot_{\alpha} x &=& bx - ay \, , \quad &z \cdot_{\alpha} y &=& cx \, ,
\quad &z \cdot_{\alpha} z &=& -ac + bz \, .
\end{matrix}
\end{equation}
To indicate the dependence on the scalars $a$, $b$, $c$, we
henceforth denote ${}^{\alpha} H_4$ by~$A_{a,b,c}$. 

It is easy to check that the center of~$A_{a,b,c}$ is trivial, i.e.,
consists of the scalar multiples of the unit.
Therefore, the two-cocycle~$\alpha$ is nondegenerate
in the sense of Section~\ref{nondegenerate} for all values of $a$, $b$,~$c$.

It follows from~\cite[Cor.~2.8]{DT2} that
$A_{a,b,c} \cong H$ when $b^2-4ac =0$ and that
$A_{a,b,c}$ is isomorphic to a quaternion algebra
when $b^2-4ac \neq 0$.
Therefore, $A_{a,b,c}$ is a simple algebra when $b^2-4ac \neq 0$.
If $b^2-4ac =0$, then $A_{a,b,c}$ is not simple
since $H$ contains a nonzero nilpotent two-sided ideal, namely the one generated by~$y$.

\subsection{The commutative algebra $\BB_{H_4}^{\alpha}$}

Let us first compute the values of the universal cocycle $\sigma$ attached to~$\alpha$,
using~\eqref{sigma-def}.

\begin{Lem}\label{sigmavalues}
(a) We have
\begin{eqnarray*}
\sigma(x,x) & = & at_x^2 t_1^{-1}, \\
\sigma(y,y) & = & \sigma(z,y) = - \sigma(y,z)
= (a t_y^2 + b t_1 t_y + c t_1^2) \, t_1^{-1},\\
\sigma(x,y) & = &  - \sigma(x,z) = (a t_x t_y - t_1t_z ) \, t_1^{-1}, \\
\sigma(y,x) & = & \sigma(z,x) =  (bt_1t_x + a t_x t_y + t_1t_z) \, t_1^{-1}, \\
\sigma(z,z) & = & -(t_z^2 + b t_x t_z + ac t_x^2) \, t_1^{-1},
\end{eqnarray*}

(b) The values of $\sigma^{-1}$ on the basis elements of~$H_4$ are linear combinations
of the values of~$\sigma$ on the same elements,
possibly divided by positive powers of~$t_1$ and of~$\sigma(x,x) = at_x^2 t_1^{-1}$.
\end{Lem}

\pf To compute the values of $\sigma$ and of~$\sigma^{-1}$, we use
\eqref{sigma-def}, \eqref{sigma^{-1}-def}, and the following values
of $\Delta^{(2)}  = (\Delta\o \id) \Delta = (\id \o \Delta\o)
\Delta$:
\begin{eqnarray*}
\Delta^{(2)}(1) & = & 1 \o 1 \o 1 \, , \\
\Delta^{(2)}(x) & = & x \o x \o x \, , \\
\Delta^{(2)}(y) & = & 1 \o 1 \o y + 1 \o y \o x + y \o x \o x \, , \\
\Delta^{(2)}(z) & = & x \o x \o z + x \o z \o 1 + z \o 1 \o 1\,  .
\end{eqnarray*}
The rest of the computation is tedious and left to the reader.
\epf

The algebra $\BB_{H_4}^{\alpha}$ is
the subalgebra of~$S(t_{H_4})_{\Theta}$ generated
by the above values of $\sigma$ and~$\sigma^{-1}$.
In order to determine it in terms of the generators of~$T(X_{H_4})$
and to obtain a presentation by generators and relations, 
we compute certain values of the universal evaluation morphism
$$\mu_{\alpha} : T(X_{H_4}) \to S(t_{H_4}) \otimes  {}^{\alpha} H_4 \, .$$
By~\eqref{universal map2} and~\eqref{H4comult}, $\mu_{\alpha}$ is
given on the generators $E$, $X$, $Y$, $Z$ by
\begin{equation}\label{H4mu}
\mu_{\alpha}(E) = t_1 \, , \quad
\mu_{\alpha}(X) = t_x\, x \, , 
\end{equation}
\begin{equation}\label{H4mu2}
\mu_{\alpha}(Y) = t_y\, x + t_1\, y \, , \quad
\mu_{\alpha}(Z) = t_z + t_x \, z \, .
\end{equation}
(In the previous formulas we have left tensor product signs
and the unit of~$H$ out.)

Set
$T = XY + YX$, $U = X(XZ + ZX)$, $V = (XZ)^2$.
Using \eqref{tablemultiplic}--\eqref{H4mu2}, we obtain the
following.

\begin{Lem}\label{muvalues}
In the algebra $S(t_{H_4})\o {}^{\alpha} H_4$ we have the following equalities:
\begin{align*}
\mu_{\alpha}(X^2) & =  a t_x^2 \, , \\
\mu_{\alpha}(Y^2) & =  a t_y^2 + b t_1 t_y + c t_1^2 \, , \\
\mu_{\alpha}(T) & =  t_x (2a t_y + b t_1) \, , \\
\mu_{\alpha}(U) & =  a t_x^2 (2t_z + b t_x) \, , \\
\mu_{\alpha}(V) & =  \mu_{\alpha}((ZX)^2)
= a t_x^2 (t_z^2 + b t_x t_z + ac t_x^2) \, , \displaybreak[1]\\
\mu_{\alpha}(4 X^2 V) & =  \mu_{\alpha}\bigl(U^2 - a^{-1}(b^2 - 4ac) \, X^6 \bigr) \, ,\\
\mu_{\alpha}(T^2 - 4 X^2 Y^2) & =  \mu_{\alpha}\bigl(a^{-1}(b^2 - 4ac)\, E^2 X^2 \bigr) \, ,\\
\mu_{\alpha}(EZ - XY) & =  t_1 t_z - a t_x t_y  \, , \\
\mu_{\alpha}(EU - X^2T) & =  2at_x^2 (t_1 t_z - a t_x t_y) \, . \\
\end{align*}
\end{Lem}

It follows from Proposition~\ref{mu-alpha-center} and Lemma~\ref{muvalues}
that the elements $E$, $X^2$, $Y^2$, $T$, $U$, $V$ of~$\UU_{H_4}^{\alpha}$ 
lie in the center~$\ZZ_{H_4}^{\alpha}$.
We henceforth identify these elements with their images
in~$\BB_{H_4}^{\alpha}$.

\begin{Thm}\label{H4universal}
(a) As a $k$-algebra,
$\BB_{H_4}^{\alpha}$ is generated by $E^{\pm 1}$, $(X^2)^{\pm 1}$, $Y^2$, $T$, $U$,
and is isomorphic to~$Z_{H_4}^{\alpha}[E^{-1},(X^2)^{-1}]$.

(b) Any relation between the generators $E^{\pm 1}$, $(X^2)^{\pm
1}$, $Y^2$, $T$, $U$ of~$\BB_{H_4}^{\alpha}$ is a consequence of the
relation
$$T^2 - 4 X^2 Y^2 = \frac{b^2- 4ac}{a} \, E^2 X^2 \, .$$
\end{Thm}

Part\,(a) of Theorem~\ref{H4universal} implies that the two-cocycle $\alpha$ is nice 
in the sense of Section~\ref{nondegenerate}
for all values of $a$, $b$,~$c$.

\pf
(a) From Lemmas~\ref{sigmavalues}\,(a) and~\ref{muvalues}
we deduce the relations
\begin{eqnarray*}
\sigma(1,1) & = & \sigma(1,x) = \sigma(x,1) = \mu_{\alpha}(E) \, , \\
\sigma(x,x) & = & \mu_{\alpha}\biggl(\frac{X^2}{E}\biggr) \, , \\
\sigma(y,y) & = & \sigma(z,y) = - \sigma(y,z) = \mu_{\alpha}\biggl(\frac{Y^2}{E}\biggr) \, ,\\
\sigma(x,y) & = &  - \sigma(x,z) = \mu_{\alpha}\biggl(\frac{X^2T - EU}{2EX^2}\biggr) \, , \\
\sigma(y,x) & = & \sigma(z,x) =  \mu_{\alpha}\biggl(\frac{X^2T + EU}{2EX^2}\biggr) \, , \\
\sigma(z,z) & = & - \mu_{\alpha}\biggl(\frac{V}{EX^2}\biggr) \, .
\end{eqnarray*}
By Lemma~\ref{sigmavalues}\,(b), the values of~$\sigma^{-1}$ are
linear combinations of the above values of~$\sigma$, possibly
divided by positive powers of $t_1 = \mu_{\alpha}(E)$ and of
$\sigma(x,x) = \mu_{\alpha}(X^2/E)$. This shows that
$B_{H_4}^{\alpha}$ is generated by $E^{\pm 1}$, $(X^2)^{\pm 1}$,
$Y^2$, $T$, $U$,~$V$. Now, by Lemma~\ref{muvalues}, $E$, $X^2$,
$Y^2$, $T$, $U$, $V$ belong to~$Z_{H_4}^{\alpha}$. Therefore,
$B_{H_4}^{\alpha}$ is obtained from $Z_{H_4}^{\alpha}$ by inverting
$E$ and~$X^2$, hence is isomorphic
to~$Z_{H_4}^{\alpha}[E^{-1},(X^2)^{-1}]$.

A computation using Lemma~\ref{muvalues} yields
$$V = \frac{U^2}{4X^2} - \frac{b^2 - 4ac}{4a} \, X^4,$$
which shows that $V$ can be expressed in terms of $(X^2)^{\pm 1}$
and~$U$.

(b) Let $B_0$ be the subalgebra of~$S(t_{H_4})_{\Theta}$ generated
by $\mu_{\alpha}(E)^{\pm 1}$, $\mu_{\alpha}(X^2)^{\pm 1}$, $\mu_{\alpha}(Y^2)$,
and $\mu_{\alpha}(U)$. Since by Part\,(a),
$$\mu_{\alpha}(E)^{\pm 1} = \sigma^{\pm 1}(1,1) \, , \quad
\mu_{\alpha}(X^2)^{\pm 1} = \mu_{\alpha}(E)^{\pm 1}\sigma^{\pm 1}(x,x) \, ,$$
$$\mu_{\alpha}(Y^2) = \mu_{\alpha}(E) \sigma(y,y) \, , \quad
\mu_{\alpha}(U) = \mu_{\alpha}(X^2) \bigl( \sigma(y,x) - \sigma(x,y) \bigr) \, ,$$
we conclude that $B_0$ is a subalgebra of~$\BB_{H_4}^{\alpha}$.
Considering the degrees in $t_1$, $t_x$, $t_y$, $t_z$ of
the polynomials $\mu_{\alpha}(E)$, $\mu_{\alpha}(X^2)$, $\mu_{\alpha}(Y^2)$,
and $\mu_{\alpha}(U/X^2)$, we easily deduce that these four polynomials are
algebraically independent.
It follows that
$$B_0 = k[E^{\pm 1}, (X^2)^{\pm 1}, Y^2, U] \, .$$

By Part\,(a), $\BB_{H_4}^{\alpha}$ is generated by~$T$ as
a~$B_0$-algebra. By Lemma~\ref{muvalues},
\begin{equation}\label{quadratic}
T^2 - 4 X^2 Y^2 = \frac{b^2- 4ac}{a} \, E^2 X^2 \in
\BB_{H_4}^{\alpha} \, .
\end{equation}
To complete the proof, it suffices to check that
$\BB_{H_4}^{\alpha}$ is a free $B_0$-module with basis $\{1, T\}$.
By~\eqref{quadratic}, the elements $1$ and $T$ generate
$\BB_{H_4}^{\alpha}$ as a $B_0$-module. Let us show that they are
linearly independent over~$B_0$. Suppose that there exists a
relation of the form $P + QT = 0$, where $P$ and $Q\in B_0$. If we
denote the degree in~$t_y$ by~$\partial_y$ and if $P, Q$ are
nonzero, we have
\begin{equation}\label{P+QT}
\partial_y (P) = \partial_y (Q ) + \partial_y (T) \, .
\end{equation}
We claim that $\partial_y (P)$ and $\partial_y (Q)$ are even
integers: indeed, of the four generators of~$B_0$ only
$\mu_{\alpha}(Y^2) = a t_y^2 + b t_1 t_y + c t_1^2$ contains~$t_y$
and its degree in~$t_y$ is~$2$. Now, $\partial_y (T) = 1$ is odd.
This contradicts~\eqref{P+QT}. Therefore, $P = Q= 0$, which shows
that $\{1, T\}$ is a basis of~$\BB_{H_4}^{\alpha}$ over~$B_0$. We
have thus proved that 
$\BB_{H_4}^{\alpha} = B_0[T]/(T^2 - 4 X^2 Y^2 - [(b^2- 4ac)/{a}] E^2 X^2)$. 
\epf

\begin{Cor}\label{H4presentation}
We have the following presentation for $\BB_{H_4}^{\alpha}$:
$$\BB_{H_4}^{\alpha}\cong
k[E^{\pm 1},(X^2)^{\pm 1},Y^2,T,U]\, / (P_{a,b,c})\, ,$$
where $P_{a,b,c}$ is the polynomial
$$ P_{a,b,c} = T^2 - 4 X^2 Y^2 - \frac{b^2- 4ac}{a} \, E^2 X^2\, .$$
\end{Cor}

It follows that, in order to specify an algebra morphism from~$\BB_{H_4}^{\alpha}$
to an extension~$K$ of~$k$,
it is enough to pick elements ${\mathbf e}$, ${\mathbf x}$, ${\mathbf y}$, ${\mathbf t}$,
${\mathbf u} \in K$ verifying ${\mathbf e} \neq 0$, ${\mathbf x} \neq 0$, and
\begin{equation*}\label{H4relation}
{\mathbf t}^2 - 4 {\mathbf x} {\mathbf y} =  \frac{b^2- 4ac}{a} \, {\mathbf e}^2 {\mathbf x} \, .
\end{equation*}

The spectrum of~$\BB_{H_4}^{\alpha}$, which is the quadric hypersurface given by
the vanishing of~$P_{a,b,c}$, is the base space of 
the universal noncommutative deformation space~$\AA_{H_4}^{\alpha}$
explicited below.

\begin{Prop}
The algebra $S(t_{H_4})_{\Theta}$ is integral over the subalgebra~$\BB_{H_4}^{\alpha}$.
\end{Prop}

\pf
It suffices to check that the generators $t^{\pm 1}_1$, $t^{\pm 1}_x$, $t_y$, $t_z$ 
of~$S(t_{H_4})_{\Theta}$ are integral over~$\BB_{H_4}^{\alpha}$. 
First, $t^{\pm 1}_1 = \sigma^{\pm 1}(1,1)$
belongs to~$\BB_{H_4}^{\alpha}$. 
For the other generators, we use the formulas in
Lemma~\ref{sigmavalues}. For instance, 
$$(t^{\pm 1}_x)^2 = a^{\mp 1} \, \sigma^{\pm 1}(1,1) \, \sigma^{\pm 1}(x,x)\, ,$$
which shows that $(t_x)^2$ and~$(t^{-1}_x)^2$ belong to~$\BB_{H_4}^{\alpha}$. 
The generator $t_y$ satisfies a quadratic equation of the form
$$a t_y^2 + b t_1 t_y + c t_1^2  - \sigma(1,1) \, \sigma(y,y) = 0$$
whose coefficients belong to~$\BB_{H_4}^{\alpha}$.
We finally use the relation
$$t_1t_z -a t_x t_y  = \sigma(1,1) \, \sigma(x,z)$$
to conclude that $t_z$ is integral over~$\BB_{H_4}^{\alpha}$.
\epf

\begin{Rem}
The elements $W = Y(YZ + ZY)$ and
$$\nabla_{\pm} = (XZY + YXZ) \pm (XYZ + YZX) $$
also belong to~$\ZZ_{H_4}^{\alpha}$. Indeed, we have
\begin{eqnarray*}
\mu_{\alpha}(W) & = & (2t_z + b t_x) (a t_y^2 + b t_1 t_y + c t_1^2)\, , \\
\mu_{\alpha}(\nabla_+) & = & t_x (2t_z + b t_x) (2a t_y + b t_1) \, , \\
\mu_{\alpha}(\nabla_-) & = & -(b^2 - 4ac)\, t_1t_x^2 \, .
\end{eqnarray*}
Moreover, the relations
$$X^2\nabla_+ = TU \quad\text{and}\quad E\nabla_- = 4 X^2 Y^2 - T^2$$
hold in~$\ZZ_{H_4}^{\alpha}$. (To prove them, it suffices to check
that both sides of each relation have equal polynomial images
under~$\mu_{\alpha}$.)
\end{Rem}

\subsection{The universal twisted algebra $\AA_{H_4}$}

Since the two-cocycle $\alpha$ is nondegenerate and nice, we may apply
Theorem~\ref{UGalois}. It states that
$$\BB_{H_4}^{\alpha} \otimes_{\ZZ_{H_4}^{\alpha}} \UU_{H_4}^{\alpha}
\quad\text{and}\quad \AA_{H_4}^{\alpha} = \BB_{H_4}^{\alpha} \otimes {}^{\sigma} H_4$$
are isomorphic as $H_4$-comodule algebras. 
Let us give a presentation of these algebras by generators and relations.

\begin{Thm}\label{AAH4}
The algebra $\AA_{H_4}^{\alpha}$ is isomorphic to
the $\BB_{H_4}^{\alpha}$-algebra generated
by two variables $\xi$, $\eta$, subject to the relations
$$\xi^2 = X^2, \quad
\eta^2 = Y^2, \quad
\xi\eta + \eta \xi = T \, .$$
\end{Thm}

\pf Let $\AA = \BB_{H_4}^{\alpha}\langle \xi, \eta\rangle /
(\xi^2 - X^2, \; \eta^2 - Y^2, \; \xi\eta + \eta \xi - T)$. 
We define an algebra morphism $f : \AA \to \AA_{H_4}^{\alpha}$ 
by $f(\xi) = X$ and $f(\eta) = Y$. It is clear that $f$ is well defined.

Let us first establish that $f$ is surjective.
The algebra $\AA_{H_4}^{\alpha}$ is generated by $E$, $X$, $Y$ and $Z$.
Now, $E$ belongs to~$\BB_{H_4}^{\alpha}$ and $X$ and $Y$ are
obviously in the image of~$f$.
From Lemma~\ref{muvalues} we deduce that
$2X^2(EZ - XY)= EU - X^2T$
holds in~$\UU_{H_4}^{\alpha}$.
Hence,
$$Z = \frac{1}{E} \, XY + \frac{EU - X^2T}{2EX^2}$$
in~$\AA_{H_4}^{\alpha}$.
Since $E^{-1}$ and $(EU - X^2T)/2EX^2$ belong to~$\BB_{H_4}^{\alpha}$,
the element $Z$ belongs to the image of~$f$.

We now prove that $f$ is injective.
We observe that $\xi^2$, $\eta^2$, $\xi\eta + \eta \xi$ belong to~$\BB_{H_4}^{\alpha}$.
Hence, $\AA$ is spanned by $1$, $\xi$, $\eta$, $\xi\eta$ as a $\BB_{H_4}^{\alpha}$-module.
Consider an arbitrary element
$\omega = \gamma_0 + \gamma_1 \xi + \gamma_2 \eta + \gamma_3 \xi \eta$
of~$\AA$, where
$\gamma_0$, $\gamma_1$, $\gamma_2$, $\gamma_3$ belong to~$\BB_{H_4}^{\alpha}$.
If $f(\omega) = 0$, then
$$\gamma_0 + \gamma_1 X + \gamma_2 Y + \gamma_3 XY = 0$$
in~$\AA_{H_4}^{\alpha}$.
Let us replace $1$, $X$, $Y$, $XY$ by their images under~$\mu_{\alpha}$.
We obtain
$$\gamma_0 \, t_1 \otimes 1 +
\gamma_1 \, (t_x \otimes x) + \gamma_2\, (t_y \otimes x + t_1
\otimes y) + \gamma_3\, (t_z \otimes 1 + t_x \otimes z) = 0 \, .$$
Since $\{1,x,y,z\}$ is a basis of~$H_4$, we obtain the system of
equations
$$\begin{cases}
\; \gamma_0 \,t_1 + \gamma_3 \, t_z & = 0 \, ,\\
\; \gamma_1 \, t_x + \gamma_2 \, t_y & =  0 \, ,\\
\; \gamma_2 \, t_1 & =  0 \, ,\\
\; \gamma_3 \, t_x & =  0
\end{cases}
$$
in~$S(t_{H_4})_{\Theta}$. 
Since $t_1$ and $t_x$ are invertible, $\gamma_i = 0$ for $i= 0, 1, 2, 3$. 
This proves the injectivity of~$f$. 
\epf

\begin{Rems}
(a) Since $A_{a,b,c}$ is simple when $b^2-4ac \neq 0$, then so is the algebra
$\KK_{H_4}^{\alpha} \otimes_{\BB_{H_4}^{\alpha}} \AA_{H_4}^{\alpha}$
by Corollary~\ref{AAsimple}.

If $b^2 - 4ac = 0$, then $P_{a,b,c} = T^2 - 4X^2 Y^2$.
Let $\theta = 2X^2 \eta - T \xi \in \AA_{H_4}^{\alpha}$, where we use the description
of~$\AA_{H_4}^{\alpha}$ given in Theorem~\ref{AAH4}.
It is easy to check that 
$$\theta^2 = \theta\xi + \xi \theta = \theta\eta + \eta \theta = 0\, .$$
Thus, $\theta$ $(\neq 0)$ generates a nilpotent two-sided ideal 
in~$\KK_{H_4}^{\alpha} \otimes_{\BB_{H_4}^{\alpha}} \AA_{H_4}^{\alpha}$,
which implies that the latter is not (semi)simple when $b^2 - 4ac = 0$.

(b) It follows from the results of Section~\ref{forms} and the computations above
that the scalar $(b^2-4ac)/a$ is an invariant for the forms of~$A_{a,b,c}$. 
More precisely, if the comodule algebra~$A_{a',b',c'}$
is a form of~$A_{a,b,c}$, then necessarily
\begin{equation}\label{discriminant}
\frac{b'{}^2- 4a'c'}{a'} = \frac{b^2- 4ac}{a} \, .
\end{equation}
It can be shown using~\cite{DT2} that, if~\eqref{discriminant} holds, 
then $A_{a',b',c'}$ is a form of~$A_{a,b,c}$.

(c) Remark (b) has an important consequence. 
Let $A_{a,b,c}$ be as above. 
Define $F=k_{0}((b^{2}-4ac)/a)$ to be the field generated by~$(b^{2}-4ac)/a$ 
over the prime field~$k_{0}$ of~$k$. 
Then there exists an $H$-comodule algebra $A_{F}$
over the field~$F$ such that $k\otimes _F A_{F} \cong A_{a,b,c}$.
Furthermore, $F$ is the unique minimal field with this property.
It follows that, if $B$ is an $H$-comodule algebra over a subfield~$K$ of 
the algebraic closure~$\overline{k}$ of~$k$  such that
$\overline{k} \otimes_K B \cong \overline{k} \otimes_k A_{a,b,c}$, 
then the universal twisted algebra $\AA_{H_{4},F}^{\alpha}$
corresponding to~$A_{F}$ specializes to~$B$. 
\end{Rems}

\appendix

\section{The map $t^{-1}$}\label{tbar-proof}

Let $C$ be a $k$-coalgebra.
Pick a copy $t_C$ of the underlying vector space of~$C$
and denote the identity map from $C$ to~$t_C$ by
$x\mapsto t_x$ ($x\in C$).
Let $S(t_C)$ be the symmetric algebra over the vector space~$t_C$
and $K_C$ the field of fractions of~$S(t_C)$.

\begin{Lem}\label{tbar}
There is a unique linear map $C \to K_C, \; x \mapsto t^{-1}_x$ such
that for all $x\in C$,
\begin{equation}\label{tbareqn}
\sum_{(x)}\, t_{x\sw1} \, t^{-1}_{x\sw2} = \sum_{(x)}\,
t^{-1}_{x\sw1}\, t_{x\sw2} = \eps(x) \, 1 \, .
\end{equation}
\end{Lem}

Here we have used the Heyneman-Sweedler sigma notation 
and $\eps$ denotes the counit of~$C$.
Observe that Equation~(\ref{tbareqn}) is equivalent to
\begin{equation}\label{tbareqn2}
t * t^{-1} = t^{-1} * t = \eps \eta \, ,
\end{equation}
where $*$ is the convolution product on $\Hom(C, K_C)$
and $\eps \eta$ is the neutral element for the convolution product
(by definition, $(\eps \eta)(x) = \eps(x) 1$ for all $x\in C$).

\pf
By~\eqref{tbareqn2} it suffices to establish that $t$ has a right and
a left inverse for the convolution product. The right and the left inverses
necessarily coincide and are unique.

Let us prove that $t$ has a right inverse for the convolution product.
First assume that $C$ is finite-dimensional.
Expressing the comultiplication $\Delta$
in a basis~$\{x_1, \ldots, x_d\}$ of~$C$, we have
\begin{equation}\label{structure constants}
\Delta(x_i) = \sum_{p,q = 0}^d\,  c_i^{p,q}\, x_p \otimes x_q \, ,
\end{equation}
where $c_i^{p,q}\in k$ are the structure constants of~$\Delta$.
When we reformulate the equations
$$\sum_{(x)}\, t_{x\sw1} \, t^{-1}_{x\sw2} = \eps(x)\, 1$$
in the basis, we obtain the matrix equation
\begin{equation}\label{tbarmatrixeqn}
M
\begin{pmatrix}
t^{-1}_{x_1}\\
\vdots\\
t^{-1}_{x_d}
\end{pmatrix}
=
\begin{pmatrix}
\varepsilon({x_1})\\
\vdots\\
\varepsilon({x_d})
\end{pmatrix}  ,
\end{equation}
where $M$ is the $d\times d$-matrix whose entries are given by
\begin{equation}\label{defmatrix-M}
M_{i,q} = \sum_{p=1}^d \, c_i^{p,q}\, t_{x_p}\in S(t_C)
\end{equation}
for all $i,q = 1, \ldots, d$.
To prove the existence and the uniqueness of the solution
of~\eqref{tbarmatrixeqn}, it suffices to
check that the determinant $\det(M)$ of~$M$ is nonzero.
Since $\det(M)$ is a polynomial in the variables~$t_{x_1}, \ldots, t_{x_d}$,
it suffices to prove that $\det(M)$ is nonzero under a suitable specialization of
the variables.

Let $\epsilon : S(t_C) \to k$ be the algebra morphism defined by
$\epsilon(t_x) = \varepsilon(x)$ for all $x\in C$.
The relation $(\varepsilon\otimes \id_C)\circ \Delta = \id_C$
satisfied by the counit becomes
$$\sum_{p,q = 0}^d\,  c_i^{p,q}\, \varepsilon(x_p) \, x_q = x_i$$
for all $i = 1, \ldots, d$. We can reinterpret this by saying that
the scalar matrix $\epsilon (M)$ whose entries $\epsilon (M)_{i,q}$ are given by
$$\epsilon (M)_{i,q} = \sum_{p=1}^d \, c_i^{p,q} \, \varepsilon(x_p)$$
satisfies
$$\epsilon (M)
\begin{pmatrix}
x_1\\
\vdots\\
x_d
\end{pmatrix}
=
\begin{pmatrix}
x_1\\
\vdots\\
x_d
\end{pmatrix}
.$$
Therefore, $\epsilon(M)$ is the identity matrix and
$\epsilon \bigl(\det(M)\bigr) =  \det \bigl(\epsilon (M)\bigr) = 1 \neq 0$.

If $C$ is not finite-dimensional, it is a direct limit of
finite-dimensional subcoalgebras~$C_{\kappa}$. By what we have just
proved, for each $\kappa$ there is a linear map $C_{\kappa} \to K_C,
\; x \mapsto t^{-1}_{\kappa, \, x}$ such that
$$\sum_{(x)}\, t_{x\sw1} t^{-1}_{\kappa, \, x\sw2} = \eps(x)\, 1$$
for all $x\in C_{\kappa}$.
Since
$$\Hom(C,K_C) = \Hom(\varinjlim_{\kappa} C_{\kappa},K_C)
\cong \varprojlim_{\kappa} \Hom(C_{\kappa},K_C) \, ,$$ 
the maps $t^{-1}_{\kappa}$ fit together to form a linear map $t^{-1}: C \to
K_C$ satisfying the required property.
One proves in an analogous way that $t$ has a left inverse for the convolution product.
\epf

We end this appendix with two elementary computations of values of~$t^{-1}$.

(a) If $g\in C$ is a grouplike element, 
i.e., satisfying $\Delta(g) = g \otimes g$ and $\eps (g) = 1$,
then 
\begin{equation}\label{grouplike-tbar}
t^{-1}_g = \frac{1}{t_g} \, .
\end{equation}

(b) If $x \in C$ is skew-primitive, i.e., if $\eps (x) = 0$ and 
$\Delta(x) = g \otimes x + x \otimes h$ for some grouplike elements~$g$,~$h$,
then
\begin{equation}\label{primitive-tbar}
t^{-1}_x = - \frac{t_x}{t_g t_h} \, .
\end{equation}

\section{The algebra $S(t_C)_{\Theta}$}\label{free-comm-coalg}

We resume the notation of Appendix~\ref{tbar-proof}, 
namely $C$ is a coalgebra, $S(t_C)$ is the symmetric algebra over 
a copy ~$t_C$ of the underlying vector space of~$C$, 
and $K_C$ the field of fractions of~$S(t_C)$.
We recall the linear maps $t: C \to S(t_C) \subset K_C$ and $t^{-1} : C \to K_C$
introduced there.

Let $S(t_C)_{\Theta}$ be the subalgebra of $K_C$ generated by 
$t(C)$ and~$t^{-1}(C)$.
The aim of this appendix is to show that $S(t_C)_{\Theta}$ is isomorphic
to the free commutative Hopf algebra on the coalgebra~$C$ introduced
by Takeuchi in~\cite[\S~11]{T}. This will imply that $S(t_C)_{\Theta}$
is obtained from~$S(t_C)$ by inverting certain ``canonical" grouplike elements.
We complete the appendix by applying the theory
to three interesting examples of coalgebras.

We first observe that, by definition of~$S(t_C)_{\Theta}$, 
the following property holds:
for any couple $(g, g^{-1})$ of linear maps $C \to R$ with values in a commutative ring~$R$ 
satisfying the relations
\begin{equation}\label{defining equation}
\sum_{(x)}\, g({x\sw1}) \, g^{-1}({x\sw2}) = \sum_{(x)}\,
g^{-1}({x\sw1}) \, g({x\sw2}) = \eps(x) \, 1
\end{equation}
for all $x\in C$, there is a unique algebra morphism 
$f : S(t_C)_{\Theta} \to R$ such that $f(t_x) = g(x)$ and 
$f(t^{-1}_x) = g^{-1}(x)$ for all $x\in C$.

For any $x\in C$, set
\begin{equation}\label{coproduct1-freecommHA}
\Delta(t_x) = \sum_{(x)}\, t_{x\sw1} \otimes t_{x\sw2} 
\in S(t_C)_{\Theta} \otimes S(t_C)_{\Theta} \, ,
\end{equation}
\begin{equation}\label{coproduct2-freecommHA}
\Delta(t^{-1}_x) = \sum_{(x)}\, t^{-1}_{x\sw2} \otimes
t^{-1}_{x\sw1} \in S(t_C)_{\Theta} \otimes S(t_C)_{\Theta} \, ,
\end{equation}
\begin{equation}\label{counit-freecommHA}
\eps(t_x) = \eps(t^{-1}_x) = \eps(x) \in k \, ,
\end{equation}
\begin{equation}\label{antipode-freecommHA}
S(t_x) = t^{-1}_x \quad\text{and}\quad  S(t^{-1}_x) = t_x \in S(t_C)_{\Theta} \, .
\end{equation}
Using \eqref{defining equation}, it is easy to check that 
\eqref{coproduct1-freecommHA}--\eqref{antipode-freecommHA}
define algebra morphisms $\Delta : S(t_C)_{\Theta} \to S(t_C)_{\Theta} \otimes S(t_C)_{\Theta}$,
$\eps: S(t_C)_{\Theta} \to k$, and 
an algebra antiautomorphism $S: S(t_C)_{\Theta} \to S(t_C)_{\Theta}$,
turning $S(t_C)_{\Theta}$ into a Hopf algebra 
with comultiplication~$\Delta$, counit~$\eps$, and antipode~$S$.

The Hopf algebra~$S(t_C)_{\Theta}$ satisfies the following universal property.

\begin{Lem}\label{universal-property}
For any commutative Hopf algebra~$H$ and any morphism of coalgebras $g : C \to H$, there is
a unique morphism of Hopf algebras $f: S(t_C)_{\Theta} \to H$
such that $f \circ t = g$.
\end{Lem}

\pf 
Since the elements $t_x$ and $t^{-1}_x$ ($x\in C$) generate
$S(t_C)_{\Theta}$ as an algebra, $f$ is uniquely defined by its
values on them. By definition, $f(t_x) = g(x)$ for all $x\in C$. 
The Hopf algebra morphism~$f$ preserves the antipodes,
which implies for all $x\in C$ that
$$f(t^{-1}_x) =  f(S(t_x))  = S(g(x))\, , $$
where the righmost~$S$ denotes the antipode of~$H$. 
To prove the existence of~$f$, it suffices by~\eqref{defining equation} to
check that
$$\sum_{(x)}\, f(t_{x\sw1}) \, f(t^{-1}_{x\sw2})
= \sum_{(x)}\, f(t^{-1}_{x\sw1}) \, f(t_{x\sw2}) = \eps(x) \, 1$$
for all $x\in C$. Now, since $g$ is a coalgebra morphism, we have
\begin{eqnarray*}
\sum_{(x)}\, f(t_{x\sw1}) \, f(t^{-1}_{x\sw2})
& = & \sum_{(x)}\, g({x\sw1}) \, S(g({x\sw2})) \\
& = & \sum_{(x)}\, {g(x)\sw1} \, S({g(x)\sw2}) \\
& = & \eps(g(x)) \, 1 = \eps(x)\, 1 \, .
\end{eqnarray*}
The relation $\sum_{(x)}\, f(t^{-1}_{x\sw1}) \, f(t_{x\sw2}) =
\eps(x)\, 1$ is proved similarly. It is now easy to check that $f$
is a morphism of Hopf algebras. 
\epf

Takeuchi's free commutative Hopf algebra generated by the coalgebra~$C$ 
(defined in~\cite[Def.~62]{T}) also satisfies the 
universal property of Lemma~\ref{universal-property}.
Therefore, it is isomorphic to~$S(t_C)_{\Theta}$.
This allows us to use the results of~\cite{T} to show that
$S(t_C)_{\Theta}$ is obtained from~$S(t_C)$
by inverting certain elements, which we now describe.

Assume first that the coalgebra $C$ is finite-dimensional.
The dual vector space $A = C^* = \Hom(C,k)$ carries a natural structure of
an algebra.
Since this algebra is finite-dimensional, we can consider the norm map
$$N : A\otimes S(t_C) \to S(t_C)$$
defined for all $\theta\in A\otimes S(t_C)$ by $N(\theta) = \det(L_{\theta})$,
where $L_{\theta}$ is
the left multiplication by~$\theta$.
On the other hand, consider the linear map
$$\iota : A \otimes S(t_C) = C^* \otimes S(t_C) \to \Hom(C,S(t_C))$$
given by $a \otimes P \mapsto (x \mapsto \langle a,x \rangle P)$,
where $a \in A$, $x\in C$, $P\in S(t_C)$, and $\langle \; , \, \rangle$
is the canonical duality between $A$ and~$C$.
The map~$\iota$ is an isomorphism since $A$ is finite-dimensional.
Using the inclusion $t\in \Hom(C,S(t_C))$
and following~\cite[\S~11]{T},
we define the element $\Theta_C \in S(t_C)$ by
\begin{equation}
\Theta_C = N\bigl( \iota^{-1}(t)\bigr) \, .
\end{equation}
By~\cite[Lemma~58]{T}, $\Theta_C$ is a grouplike element of~$S(t_C)$.

We shall see in Example~\ref{exempleS(k^G)} below that
$\Theta_C$ is a generalization of Dede\-kind's group determinant.
For this reason, to denote this element
we have used the same symbol~$\Theta$ as Dedekind.
Note also that $\Theta_C$ is the determinant of the matrix~$M$
defined by~\eqref{defmatrix-M}.

In the general case, we have the following.

\begin{Prop}\label{Takeuchi}
Let $C$ be a coalgebra and $(C_{\kappa})_{\kappa}$ be a family of finite-dimensional
coalgebras such that $\sum_{\kappa}\, C_{\kappa}$ contains the coradical~$C_0$ of~$C$.
Then $S(t_C)_{\Theta}$ is the smallest subalgebra of the fraction field~$K_C$
containing~$S(t_C)$ and all inverses~$\Theta^{-1}_{C_{\kappa}}$:
$$S(t_C)_{\Theta} =  S(t_C)\biggl[\biggl(\frac{1}{\Theta_{C_{\kappa}}}\biggr)_{\kappa} \, \biggr] \, .$$
Moreover, if the coradical $C_0$ is finite-dimensionsal, then
$$S(t_C)_{\Theta} =  S(t_C)\biggl[\frac{1}{\Theta_{C_0}} \biggr] \, .$$
\end{Prop}

\pf
This is a consequence of \cite[Theorem~61 and Corollary~64]{T} together
with our identification of~$S(t_C)_{\Theta}$ with
Takeuchi's free commutative Hopf algebra on~$C$.
\epf

As a consequence of Proposition~\ref{Takeuchi}, the spectrum of the commutative
algebra $S(t_C)_{\Theta}$ is obtained from the affine space~$C$ by cutting out
the zero-loci of the polynomials $\Theta_{C_{\kappa}}$.
We can reinterpret this as follows: for any field~$K$,
a linear map $g : C\to K$ extends to an algebra morphism 
$S(t_C)_{\Theta} \nolinebreak \to \nolinebreak K$
if and only if its extension $f : S(t_C) \to K$ as an algebra morphism
satisfies $f(\Theta_{\kappa}) \neq 0$ for all~$\kappa$.

We end this appendix with three examples.

\begin{Expl}\label{exempleM_n(k)}
Let $C$ be the coalgebra dual to the algebra $M_n(k)$ of $n\times n$-matrices
with entries in~$k$. 
As a vector space, $C$ has a basis $\{X_{i,j}\}_{i,j = 1, \ldots, n}$
with comultiplication~$\Delta$ and counit~$\eps$ given by
$$\Delta(X_{i,j})  = \sum_{k=0}^n\, X_{i,k} \o X_{k,j}
\quad \text{and} \quad \eps(X_{i,j}) = \delta_{i,j}\;\;
\text{(Kronecker symbol)}$$ 
for all $i,j = 1, \ldots, n$. 
Set $t_{i,j}  = t_{X_{i,j}}$. 
Then $\Theta_C \in S(t_C)$ is equal to the determinant of
the ``generic" matrix $\bigl(t_{i,j}\bigr)_{i,j}\,$.
\end{Expl}

\begin{Expl}\label{exempleS(k[G])}
Let $C$ be a coalgebra which has a basis~$G$ consisting of grouplike
elements (for instance, the coalgebra underlying the Hopf
algebra of a group~$G$ considered in Example~\ref{graded algebra}). 
The symmetric algebra $S(t_C)$ is the polynomial algebra
in the indeterminates~$(t_g)_{g\in G}$. 
By~\eqref{grouplike-tbar}, $t^{-1}_g = 1/t_g$ for all $g\in G$. 
Therefore, $S(t_C)_{\Theta}$ is the algebra of Laurent polynomials
$$S(t_C)_{\Theta} = k[t_g,t_g^{-1}\, |\;  g\in G] \, .$$
Observe that $\Theta_C = \prod_{g\in G}\, t_g$ if $G$ is finite.
\end{Expl}

\begin{Expl}\label{exempleS(k^G)}
Let $G$ be a finite group and $C = k^G$ the coalgebra underlying the Hopf
algebra of $k$-valued functions on~$G$
considered in Example~\ref{algebra with G-action}.
(Over an algebraically closed field, the coalgebra~$k^G$
is a product of matrix coalgebras of type~\ref{exempleM_n(k)}.)
In this case $\Theta_C$ is equal to {\em Dedekind's group determinant}~$\Theta_G$,
which is the determinant of the matrix
$$\Bigl(t_{gh^{-1}} \Bigr)_{g,h \in G}\, .$$
The determinant~$\Theta_G$ was factored into a product of irreducible polynomials
by Dedekind when $G$ is abelian and by Frobenius for an arbitrary group.
The factors can be expressed
in terms of the irreducible characters of~$G$, see~\cite[Chap.~12, Part~B]{vdW}.

If $G$ is an abelian group of exponent~$N$, and $k$ is a field
containing a primitive $N$-th root of unity and having
a characteristic not dividing~$N$,
then there is an isomorphism of Hopf algebras $C = k^G \cong k[G]$
provided by the discrete Fourier transform.
Under these hypotheses, the algebra $S(t_C)_{\Theta}$ is
a ring of Laurent polynomials in view of Example~\ref{exempleS(k[G])}.

If $G = \langle g \; |\, g^2 = e\rangle$
is the cyclic group of order two and $C = k^G$, then
$$\Theta_C = 
\left|
\begin{array}{cc}
t_e & t_g \\
t_g & t_e
\end{array}
\right| 
= t_e^2 - t_g^2 = (t_e - t_g)(t_e + t_g) \, .$$
If $k$ is of characteristic two,
then $S(t_C)_{\Theta} = k[t_e, t_g, (t_e + t_g)^{-2}]$.
One can check that this is not an algebra of Laurent polynomials.
\end{Expl}


\begin{thebibliography}{99}

\bibitem{Abe} E.~Abe, 
\textit{Hopf algebras},
Cambridge Tracts in Mathematics,~74,
Cambridge University Press, Cambridge-New York,~1980.

\bibitem{AEGN} E.~Aljadeff, P.~Etingof, S.~Gelaki, D.~Nykshych,
\textit{On twisting of finite-dimensional Hopf algebras},
J.~Algebra 256 (2002), 484--501.


\bibitem{AHN1} E.~Aljadeff, D.~Haile, M.~Natapov,
\textit{Graded identities of matrix algebras and the universal graded algebra},
preprint~2007.

\bibitem{AN} E.~Aljadeff, M.~Natapov,
\textit{On the universal $G$-graded central simple algebra},
pre\-print~2007.


\bibitem{Az} G.~Azumaya,
\textit{On maximally central algebras},
Nagoya Math.\ J.~2 (1951), 119--150.

\bibitem{BL} Y. A. Bahturin, V. Linchenko,
\textit{Identities of algebras with actions of Hopf algebras},
J.~Algebra 202 (1998), 634--654.

\bibitem{BZ} Y. A. Bahturin, M.~Zaicev,
\textit{Identities of graded algebras},
J.~Algebra 205 (1998), 1--12.

\bibitem{BC} J.~Bichon, G.~Carnovale, 
\textit{Lazy cohomology: an analogue of the Schur multiplier for arbitrary Hopf algebras},
J.~Pure Appl.\ Algebra 204 (2006), 627--665.

\bibitem{BCM} R. J. Blattner, M. Cohen, S. Montgomery, 
\textit{Crossed products and inner actions of Hopf algebras},
Trans.\ Amer.\ Math.\ Soc.\ 298 (1986), 671--711. 

\bibitem{BM} R. J. Blattner, S. Montgomery, 
\textit{A duality theorem for Hopf module algebras},
J.~Algebra 95 (1985), 153--172.


\bibitem{BDZ} T.~Brzezi\'nski, L.~D\c abrowski, B.~Zieli\'nski, 
\textit{Hopf fibration and monopole connection over the contact quantum spheres},
J.~Geom.\ Phys.\ 50 (2004), 345--359. 

\bibitem{BM1}  T.~Brzezi\'nski, S. Majid,
\textit{Quantum gauge group theory on quantum spaces},
Comm.\ Math.\ Phys.\ 157 (1993), 591--638; Erratum: 167 (1995), 235.

\bibitem{BM2} T.~Brzezi\'nski, S. Majid,
\textit{Line bundles on quantum spheres},
in Particles, fields, and gravitation (\L \'od\'z 1998), AIP Conf.\ Proc., 453,
Amer.\ Inst.\ Phys., Woodbury, New York, 1998, 3--8.

\bibitem{Doi} Y. Doi, 
\textit{Equivalent crossed products for a Hopf algebra},
Comm.\ Algebra 17 (1989), 3053--3085.

\bibitem{DT1} Y.~Doi, M.~Takeuchi, 
\textit{Cleft comodule algebras for a bialgebra},
Comm.\ Algebra 14 (1986), 801--817.

\bibitem{DT2} Y. Doi, M. Takeuchi, 
\textit{Quaternion algebras and Hopf crossed products},
Comm.\ Algebra 23 (1995), 3291--3325.

\bibitem{Eis} D. Eisenbud, 
\textit{Commutative algebra. With a view toward algebraic geometry}, 
Grad.\ Texts in Math.~150, Springer-Verlag, New York,~1995.

\bibitem{GMS} S.~Garibaldi, A.~Merkurjev, J.-P.\ Serre, 
\textit{Cohomological invariants in Galois coho\-mol\-ogy}, 
Univ.\ Lecture Ser.~28, Amer.\ Math.\ Soc., Providence, RI, 2003.

\bibitem{Ha1} P. M. Hajac,
\textit{Strong connections on quantum principal bundles},
Comm.\ Math.\ Phys.\ 182 (1996), 579--617.

\bibitem{HM} P. M. Hajac, S. Majid,
\textit{Projective module description of the $q$-monopole},
Comm.\ Math.\ Phys.\ 206 (1999), 247--264.

\bibitem{K} C.~Kassel, 
\textit{Quantum principal bundles up to homotopy equivalence},
The Legacy of Niels Henrik Abel, The Abel Bicentennial, Oslo, 2002,
O.~A.~Laudal, R.~Piene (eds.), Springer-Verlag 2004, 737--748
(see also arXiv:math.QA/0507221).

\bibitem{KS} C.~Kassel, H.-J.\ Schneider,
\textit{Homotopy theory of Hopf Galois extensions},
Ann.\ Inst.\ Fourier (Grenoble) 55 (2005), 2521--2550.

\bibitem{LPR} G.~Landi, C.~Pagani, C.~Reina,
\textit{A Hopf bundle over a quantum four-sphere from the symplectic group},
Comm.\ Math.\ Phys.\ 263 (2006), 65--88. 

\bibitem{La} S.~Lang, 
\textit{Algebra},
Addison-Wesley Publishing Co., Inc., Reading, Mass.,~1965.

\bibitem{Ma1} A. Masuoka,
\textit{Cleft extensions for a Hopf algebra generated by a nearly primitive element},
Comm.\ Algebra 22 (1994), 4537--4559.

\bibitem{M2} S. Montgomery, 
\textit{Hopf algebras and their actions on rings},
CBMS Conf.\ Series in Math., vol.~82, Amer.\ Math.\ Soc., Providence, RI, 1993.

\bibitem{Row} L.~Rowen, 
\textit{Polynomial identities in ring theory},
Academic Press, New York,~1980.

\bibitem{S} H.-J.\ Schneider, 
\textit{Principal homogeneous spaces for arbitrary Hopf algebras}, 
Israel J.\ Math.\ 72 (1990), 167--195.

\bibitem{Sw} M.~Sweedler, 
\textit{Hopf algebras}, 
W. A. Benjamin, Inc., New York,~1969.

\bibitem{T} M.~Takeuchi, 
\textit{Free Hopf algebras generated by coalgebras},
J.~Math.\ Soc.\ Japan 23 (1971), 561--582.

\bibitem{vdW} B.-L.\ van der Waerden,
\textit{A history of algebra. From al-Khw\=arizm\=\i\ to Emmy Noether},
Springer-Verlag, Berlin,~1985.



\end{thebibliography}
\end{document}